\documentclass[12pt]{amsart}

\usepackage[pagebackref,colorlinks=true]{hyperref}  

\usepackage{mathtools}

\usepackage{amsfonts}
\usepackage{amsmath,amsthm,amssymb,amscd,enumerate,eucal,url}

\usepackage{eucal,url,amssymb,verbatim,
	enumerate,amscd,
}

\allowdisplaybreaks[1]

\numberwithin{equation}{section}

\newtheorem{thrm}{Theorem}[section]
\newtheorem{lemma}[thrm]{Lemma}
\newtheorem{prop}[thrm]{Proposition}
\newtheorem{cor}[thrm]{Corollary}
\newtheorem{dfn}[thrm]{Definition}

\newtheorem{rmrk}[thrm]{Remark}

\newtheorem*{conv*}{Conventions}

\usepackage[margin=1in]{geometry}
\newcommand{\R}{\mathbb{R}}
\newcommand{\lc}{\langle}
\newcommand{\rc}{\rangle}
\newcommand{\br}[1]{\overline{#1}}
\newcommand{\dt}{\,.\,}
\newcommand{\spc}{\ \,}

\newcommand{\ind}[1]{\text{\scalebox{0.85}{\ensuremath #1}}}
\newcommand{\pow}[1]{\text{\scalebox{1.15}{\ensuremath #1}}}

\begin{document}
	
	\setlength{\abovedisplayskip}{8pt} 
	\setlength{\belowdisplayskip}{8pt}
	
	\begin{abstract}  This paper examines 8-dimensional Riemannian manifolds whose structure group reduces to ${SO(4)}_{ir}\subset GL(8,\mathbb R)$, the image of an irreducible representation of $SO(4)$ on $\mathbb R^8$.  We demonstrate that such a reduction can be described by an almost quaternion-Hermitian structure and a special rank-4 tensor field, which we call a cubic discriminant. This tensor field is pointwise linearly equivalent to the formula for the discriminant of a cubic polynomial. 
		We show that the only non-flat, integrable examples of these structures are the quaternion-K\"ahler symmetric spaces  $G_2\big/SO(4)$ and $G_{2(2)}\big/SO(4)$. We also present a new curvature-based characterization for the Riemannian metrics on these spaces.
	\end{abstract}
	
	
	
	\title[A geometry of cubic discriminants in 8 dimensions]
	{A geometry of cubic discriminants in 8 dimensions}
	
	\date{\today}
	
	\author{Elitza Hristova}
	\address[Elitza Hristova]{
		Institute of Mathematics and Informatics,
		Bulgarian Academy of Sciences,
		Acad. G. Bonchev Str., Block 8,
		1113 Sofia, Bulgaria}
	\email{e.hristova@math.bas.bg}

	\author{Ivan Minchev}
	\address[Ivan Minchev]{ Faculty of Mathematics and Informatics, University
		of Sofia, blvd. James Bourchier 5, 1164 Sofia, Bulgaria}
	\email{minchev@fmi.uni-sofia.bg}
	
	\maketitle

	\setcounter{tocdepth}{1} \tableofcontents

	\section{Introduction}
	The Lie group $SO(4)$ has two non-isomorphic, irreducible 8-dimensional representations on $\mathbb R^8$. Because these representations are connected by an automorphism of  $SO(4)$, they generate the same image subgroup, ${SO(4)}_{ir}\subset GL(8,\mathbb R),$  up to conjugation.
	In the paper we examine 8-dimensional manifolds where the structure group of the bundle of linear frames can be reduced to ${SO(4)}_{ir}$. We demonstrate that any such reduction is defined by two components: an almost quaternion-Hermitian structure on the manifold and a unique rank-4 tensor field. This tensor field, which we call a cubic discriminant, is, in a specific sense, pointwise linearly equivalent to the formula for the discriminant of a cubic polynomial. Our approach is similar to the study of 5-dimensional manifolds where the structure group is reduced to the irreducible $SO(3)$ (cf. \cite{BN, ABF, CF}). In that case, the geometry is described by a Riemannian metric and a special rank-3 tensor field.

	We begin with the space of complex cubic polynomials, $W=\mathbb C^4,$ where each polynomial is of the form $p(z)=az^3+bz^2+cz+d$. This space is equipped with the standard complex representation of the group $Sp(2)$. The discriminant of $p(z)$, given by the formula:
	\begin{equation}\label{cubic_dis}
		Dis(a,b,c,d)=18abcd-27a^2d^2-4ac^3-4b^3d+b^2c^2
	\end{equation}
	defines a specific element, $\hat S,$ within the space of all totally symmetric covariant 4-tensors on $W$ ($S^4\,W^\ast$). 
	If we consider the standard complex representation of $Sp(1)$ on $H=\mathbb C^2$, then the tensor product $V=W\otimes H$ inherits a representation of the product group $Sp(2)Sp(1)$. Interestingly, the tensor $\hat S\in S^4\,W^\ast$ (defined by discriminant \eqref{cubic_dis}) can also be interpreted as a hyper-K\"ahler curvature type tensor on $V.$ This correspondence allows us to define a unique $Sp(2)Sp(1)$-orbit, $\mathcal C,$ in the space of such tensors on $V.$ In this paper, we provide a coordinate-free algebraic description of $\mathcal C$ (Theorem~\ref{sp2-characterization}) and prove that it is a 7-dimensional homogeneous space, specifically:
	\begin{equation}\label{sp2sp1-so4_0}
		\mathcal C\cong\frac{Sp(2)Sp(1)}{SO(4)_{\text{ir}}}.
	\end{equation}
	
	A cubic discriminant on an almost quaternion-Hermitian 8-manifold $M$ is, by definition, a section of a certain locally trivial fiber bundle $\mathcal C(M)\rightarrow M$ whose fibers are diffeomorphic to \eqref{sp2sp1-so4_0}. There is a one-to-one correspondence between these cubic discriminants on $M$ and the possible reductions of the frame bundle of $TM$ from the structure group $Sp(2)Sp(1)$ to the subgroup $SO(4)_{ir}$. The presence of a cubic discriminant on an 8-dimensional almost quaternion-Hermitian manifold $M$ allows us to treat its tangent vectors, in a sense, as cubic polynomials. To clarify, let us first assume the structure group of the tangent bundle can be reduced from $Sp(2)Sp(1)$ to $Sp(2)$. This means that the manifold $M$ has a global triple of compatible almost complex structures. While this assumption might not hold globally, it is always satisfied locally. In this case, choosing a cubic discriminant on $M$ is equivalent to identifying the complex tangent bundle $TM$ ($TM$ is considered complex with respect to one of the compatible almost complex structures) with  $\text{S}^3\,\Delta^\ast$ for a certain complex vector bundle $\Delta\rightarrow M$ with 2-dimensional fibers. This means that each tangent vector is a rank-3 symmetric polynomial on the respective (projectivized) fiber of $\Delta$. In order to be able to upgrade this interpretation  to the general case (of a structure group $Sp(2)Sp(1)$), one needs to involve the concept of a twistor space, cf. \cite{Sal,LeB}. This clearly is a geometrically very rich setup for various considerations. For example, we believe this can be used to introduce a meaningful first-order Dirac-like operator on the sections of $\Delta$ by replacing the Clifford multiplication between tangent vectors and spinors in the classical setting (see for example \cite{Fr}) with the obvious action of $\text{S}^3\,\Delta^\ast$ on $\Delta$. However, in this paper, we shall not explore this line of thought any further but rather focus on some more fundamental issues concerning the geometry of cubic discriminants.
	
	Section~\ref{sec_can_con} establishes the existence of a unique canonical linear connection on $TM$ for any given cubic discriminant on a quaternion-K\"ahler manifold. This connection has torsion that splits equivariantly into four irreducible components. We present two natural examples of homogeneous cubic discriminants on quaternion-K\"ahle manifolds: the compact Wolf space $G_2/SO(4)$ \cite{W}; and the non-compact symmetric space $G_{2(2)}/SO(4)$ \cite{W,Ale,Cor}, where $G_2$ and $G_{2(2)}$ are respectively the compact  and the split forms of the complex $G^c_2$. We prove that these two are locally the only integrable non-flat cubic discriminants (Theorem~\ref{integrable_cubics}).  
	
	A quaternion-K\"ahler manifold's Riemannian curvature tensor, $R$, has a well-known decomposition: $R=R'+\frac{Scal}{64}R_0$ (see Section~\ref{almostQH} below for more details). Here is a breakdown of the components:
	\begin{itemize}
		\item $R'$ is a hyper-K\"ahler curvature type tensor.
		\item $Scal$ is the manifold's constant scalar curvature.
		\item $R_0$ is a specific tensor field that is defined by the curvature of the quaternionic projective space, $\mathbb HP^2$.
	\end{itemize}
	If $R'$ is zero, the manifold is locally isometric to one of three spaces: the quaternion projective space ($\mathbb HP^2$), the quaternion hyperbolic space ($\mathbb HH^2$), or the flat space ($\mathbb H^2$). In the paper, we demonstrate (Theorem~\ref{main_thrm_third}) that, rather than assuming $R'$ vanishes, if we suppose that $R'$ is pointwise proportional to a cubic discriminant, then $M$ is locally isometric to either $G_2/SO(4)$ or $G_{2(2)}/SO(4)$. This provides a novel curvature characterization for the Riemannian metrics on these two exceptional symmetric spaces.

	\begin{conv*}
		In the paper we use the following general conventions:
		\begin{enumerate}[a)]
			\item The small  Latin indices $s,t,\dots$ are usually running from 1 to 3 (when nothing else specified).
			\item The conventions related to the small Greek indices $\alpha,\beta,\gamma,\dots$ are explained in section~\ref{index_conventions}
			\item We will often, but not always, use without comment the convention of summing over repeated indices.
			\item If $V$ is a vector space (or a bundle over a manifold),  then $V^\ast$ is its dual. $\text S^k\,V$ and $\Lambda^k\,V$ denote respectively the $k^{\text{th}}$ totally symmetric and skew-symmetric tensor power of $V$. $End(V)$ is the space (or bundle) of endomorphisms of $V,$ $End(V)\cong V^\ast\otimes V$.
		\end{enumerate}
	\end{conv*}

	\textbf{Acknowledgments}  
	The research of E.H. is partially supported by the Bulgarian Ministry of Education
	and Science, Scientific Programme ”Enhancing the Research Capacity in Mathematical Sciences (PIKOM)”, No. D01-88/20.06.2025 and by COST (European Cooperation in Science and Technology) via COST Action 21109 CaLISTA.
	The research of I.M.  is partially supported by Contract KP-06-H72-1/05.12.2023 with the National Science Fund of Bulgaria and  by Contract 80-10-61/27.05.2025  with the Sofia University "St.Kl.Ohridski" and  the National Science Fund of Bulgaria.

	\section{Preliminaries}\label{preliminaries}
	
	A hyper-Hermitian structure  on a real $4n$-dimensional (for the most part $n$ will be $2$) vector space $V$ with a positive definite inner product $g$ is defined by a triple of endomorphisms $J_1,J_2,J_3\in \text{End}(V)$ that satisfy:
	\begin{equation}\label{quat_iden}
		J_1\cdot J_2=-J_2\cdot J_1=J_3;\qquad J_s\cdot J_s=-\text{Id};\qquad\text{and}\qquad g(J_sx,J_sy)=g(x,y),
	\end{equation}
	for all $x,y\in V$ and $s=1,2,3$. The group of all linear transformations preserving a fixed hyper-Hermitian structure is known to be isomorphic to the compact symplectic group $Sp(n)$. Alternatively, on a complex $2n$-space $W$ with a prescribed non-degenerate 2-from $\pi$, a hyper-Hermitian structure can be defined by
	choosing an anti-linear transformation $j$ of $W$ that satisfies:
	\begin{equation}\label{q-str}
		j^2=-\text{Id};\qquad \pi(jx,jy)=\overline{\pi (x,y)};\qquad\text{and}
		\qquad \pi(x,jx)>0\ \ \text{if}\ \ x\ne 0.
	\end{equation}
	
	To clarify the relationship between \eqref{quat_iden} and \eqref{q-str}, suppose that we are given the latter, i.e., a triple $(W,\pi, j)$ with the above properties. Then $V$ is equal to the underlying real vector space of $W$. We define the complex structure $J_1 :V \rightarrow V$ as $J_1(v) = iv$ for all $v \in V$, where $i$ is the imaginary unit of $W$. Notice that, if we consider the complexification of $V$, namely the space $V^{\mathbb{C}} = V \otimes_{\mathbb{R}}\mathbb{C}$, and we extend $J_1$ by linearity to $V^{\mathbb{C}}$ via the formula $J_1(v \otimes z) = J_1(v) \otimes z$ for all $v \in V$ and $z \in \mathbb{C}$, then $W$ is naturally isomorphic to the $i$-eigenspace of $J_1$. In addition, if $\overline{W}$ denotes the complex conjugate vector space of $W$, then $\overline{W}$ is isomorphic to the $(-i)$-eigenspace of $J_1$ and we have the splitting 
	
	\begin{equation}\label{splitVc}
		V^{\mathbb C} \cong W\oplus\overline{W}.
	\end{equation} 
	
	
	Next, let $c : W \rightarrow \overline{W}$ denote the complex conjugation map. We define $J_2 = j \circ c$. Then $J_2$ is a complex linear map from $W$ to $\overline{W}$. We take the corresponding conjugated map $\bar{J_2}: \overline{W} \rightarrow W$ and obtain the complex linear map:
	\[
	J_2 \oplus \bar{J_2} : W\oplus \overline{W} \rightarrow \overline{W} \oplus W. 
	\]
	With a slight abuse of notation we denote $J_2 = J_2\oplus \bar{J_2}$. Then, we define $J_3 = J_1J_2$.
	
	Next, we define a symmetric bilinear form $g: V^{\mathbb{C}} \times V^{\mathbb{C}} \rightarrow \mathbb{C}$ by the formula
	\begin{equation}\label{g-pi}
		g(x,y)\overset{def}{\ =\ }\pi(x,J_2y).
	\end{equation}
	It is easily seen that $J_1, J_2,J_3$ and $g$ satisfy the hyper-Hermitian identities \eqref{quat_iden} and thus the stabilizer group of the pair $(j,\pi)$ in the general linear group $GL_{\mathbb C}(W)$ is precisely $Sp(n)$. 
	
	The linear span 
	\begin{equation}\label{spanJ}
		\mathcal Q=\text{span}_{\mathbb R}\{J_1,J_2,J_3\}
	\end{equation}
	is a three dimensional subspace of $End(V)\footnotemark$\footnotetext{$End(V)$ denotes the space of all endomorphisms of $V$} that has a natural inner product and an orientation with respect to which $J_1,J_2,J_3$ is a positively oriented orthonormal basis. In fact, one can easily check that a triple of elements of $\mathcal Q$ satisfies the identities \eqref{quat_iden} if and only if it constitutes a positively oriented orthonormal basis of $\mathcal Q$. We have that: 
	\begin{equation}\label{def_spnsp1}
		Sp(n)Sp(1)=\Big\{A\in SO(V)\ \Big|\ A\mathcal QA^{-1}\subset \mathcal Q\Big\},
	\end{equation}
	where 
	$$Sp(1)=\Big\{x_0\,\text{Id}+\sum_{s=1}^3x_sJ_s\ \Big|\ \sum_{s=0}^3x_s^2=1\}.$$

	\subsection{Conventions concerning the use of complex indices in the paper}\label{index_conventions}  We will be using a fixed hyper-Hermitian structure $(W,j,\pi)$, defined as in \eqref{q-str}.  From this, we construct the vector space $V,$ the metric $g$ and the complex structures $J_s$, $s=1,2,3,$ as outlined in \eqref{splitVc} and \eqref{g-pi}.
	
	If $\{e_\alpha\}$ is any basis of $W$,  then \mbox{$\{e_{\bar\alpha} \overset{def}{\ =\ }\overline {e_\alpha}\}$} is a basis of $\overline W$. We shall denote by $\{e^\alpha\}$ and  $\{e^{\bar\alpha}\}$ the respective dual bases for $W^\ast$ and $\overline{W^{\ast}}$ (obviously, we are distinguishing between upper and lower indices). Setting   $\pi_{\alpha\beta}=\pi(e_\alpha,e_\beta)$ and $g_{\alpha\bar\beta}=g(e_\alpha,e_{\bar\beta} )$, we let $g^{\alpha\bar\beta}$ be the matrix inverse of $g_{\alpha\bar\beta}$.  
	In the sequel, we shall use without comment the usual summation convention over repeating indices.  The small Greek indices $\alpha,\beta,\gamma,\dots$will be running from $1$ to $2n$ (for the most part of the paper $n$ will be $2$). Any array of complex numbers indexed by lower and upper Greek indices (with and without bars) that appears in the text will be considered as a tensor on $V^{\mathbb C}$ (with respect to the fixed basis), e.g., $A_{\alpha\dt\dt}^{\spc\beta\bar\gamma}$ corresponds to the tensor 
	\begin{equation*}
		A_{\alpha\dt\dt}^{\spc\beta\bar\gamma}\, e^\alpha\otimes e_{\beta}\otimes e_{\bar\gamma}.
	\end{equation*}     
	Clearly, the vertical as well as the horizontal position of an index carries information about the tensor.
	The matrices $g_{\alpha\bar\beta}$ and $g^{\alpha\bar\beta}$ will be used to lower and raise indices in the usual way, e.g.,  
	\begin{equation}\label{index-r-l}
		x_\alpha=g_{\alpha\bar\beta}x^{\bar\beta},\qquad  A^{\alpha}_{. \, \beta}=g^{\alpha\bar\gamma}A_{\bar\gamma\beta}.
	\end{equation}
	Another important rule that we follow here is that each tensor that we consider must be real, i.e., a tensor on the real component $V$ of $V^{\mathbb C}=W\oplus\overline{W}
	$; this means that whenever an array like $A_{\alpha\ \ }^{\spc\beta\bar\gamma}$ appears,  the array $A_{\bar\alpha\ \ }^{\spc\bar\beta\gamma}$ is assumed to be defined, by default, via complex conjugation,
	\begin{equation}\label{real_t}
		A_{\bar\alpha\ \ }^{\spc\bar\beta\gamma}= \overline{A_{\alpha\ \ }^{\spc\beta\bar\gamma}}.
	\end{equation}
	In fact, we shall sometimes omit the explicit tensoring with $\mathbb C$  from the notation and let it be understood from the context whether $V$ denotes a real vector space or its complexification. 
	
	Following these conventions, the complex structure $J_2$ on $V$ is represented by an array $(J_2)^{\alpha}_{.\, \bar{\beta}}.$ By \eqref{g-pi}, we have 
	\[
	\pi^{\alpha}_{.\, \bar\beta} = -(J_2)^{\alpha}_{.\, \bar\beta}
	\]
	and therefore,
	\begin{equation*}
		\pi^{\alpha}_{.\, \bar\sigma}\pi^{\bar\sigma}_{.\, \beta}=
		\begin{cases}
			-1,\quad  \text{if}\ \alpha=\beta\\
			\ \ 0,\quad \text{otherwise}.
		\end{cases}
	\end{equation*}

	We shall say that a basis $\{e_\alpha\}$ of $W$ is $Sp(n)$-adapted if
	\begin{equation}\label{constants}
		e_{n+a}=j(e_a),\quad a=1,\dots n, \qquad \text{and}\qquad 
		\pi = \sum_{a=1}^n e^a\wedge e^{n+a}.
	\end{equation}
	\noindent For any $Sp(n)$-adapted basis, we have
	\begin{equation}\label{g_ab}
		g_{\alpha\bar\beta}\ =\ g( e_\alpha,e_{\bar\beta})\ =\ \pi\Big(e_\alpha,j(e_{\beta})\Big)\ =\  \begin{cases}1, & \mbox{if } \alpha=\beta\\0, & \mbox{if } \alpha\ne\beta.
		\end{cases}
	\end{equation}

	\subsection{A few facts about hyper-Hermitian structures}
	
	We introduce a complex antilinear transformation $\mathfrak j$ on the tensor algebra of $V^{\mathbb{C}}$ as follows: 
	$\mathfrak{j}(v) = -J_2(\bar v)$ for all $v \in V^{\mathbb{C}}$. Then, for an arbitrary element $\sum_a v_{a1}\otimes \cdots \otimes v_{ar} \in (V^{\mathbb{C}})^{\otimes r}$ we set:
	\[
	\mathfrak{j}(\sum_a v_{a1}\otimes \cdots \otimes v_{ar}) = \sum_i\mathfrak{j}(v_{a1}\otimes \cdots \otimes v_{ar}) = \sum_a \mathfrak{j}(v_{a1})\otimes \cdots \otimes \mathfrak{j}(v_{ar}).
	\]
	In coordinates, $\mathfrak j$ takes any tensor $T$ with components $T_{\alpha_1\dots\alpha_k\bar\beta_1\dots\bar\beta_l\dots}$ to a tensor $\mathfrak jT$ of the same type, whose components are given by 
	\begin{equation}\label{map_j}
		(\mathfrak jT)_{\alpha_1\dots\alpha_k\bar\beta_1\dots\bar\beta_l\dots}=\sum_{\bar\sigma_1\dots\bar\sigma_k\tau_1\dots\tau_l\dots}\pi^{\bar\sigma_1}_{\  \alpha_1}\dots\pi^{\bar\sigma_k}_{\  \alpha_k}\,\pi^{\tau_1}_{\  \bar\beta_1}\dots\pi^{\tau_l}_{\  \bar\beta_l}\dots T_{\bar\sigma_1\dots\bar\sigma_k\tau_1\dots\tau_l\dots}.
	\end{equation}
	Clearly, $\mathfrak j$ is a $Sp(n)$-equivaraint map and the restriction of $\mathfrak j$ to $W$ coincides with $j$ \big(cf.~\eqref{q-str}\big). Using this, we have the following useful characterization for the Lie algebra $sp(n)$ of the compact symplectic Lie group $Sp(n)$. 
	\begin{lemma}\label{sp(2)}
		Each of the following two vector spaces (for the description of which we use any fixed basis $\{e_\alpha\}$ of $W$ and the conventions introduced in Section~\ref{index_conventions}) is isomorphic to $sp(n)$:
		
		\vspace{0.3cm}\par (1) The subspace of  $\text{S}^2\,{W^*}$ given by\footnotemark\footnotetext{$\text{S}^2\,{W^*}$ is the space of all (complex)  symmetric covariant 2-tensors on $W$}
		\begin{equation}\label{Spn_X}
			\Big\{X_{\alpha\beta}\ \Big|\ X_{\alpha\beta}=X_{\beta\alpha},\ (\mathfrak jX)_{\alpha\beta}=X_{\alpha\beta}\Big\}
		\end{equation}
		\par (2) The subspace of $End(V)$ given by\footnotemark\footnotetext{$End(V)$ is the space of all endomorphisms of $V,$ $End(V)=V^\ast\otimes V.$}
		\begin{equation}\label{Spn_A}
			\Big\{A_{\alpha\bar\beta}\ \Big|\ A_{\alpha\bar\beta}=-A_{\bar\beta\alpha},\ (\mathfrak jA)_{\alpha\bar\beta}=A_{\alpha\bar\beta}\Big\}
		\end{equation}
		
		\vspace{0.5cm}\par The equivalence between  (1) and (2) is explicitly given by  $A^\alpha_{\dt\beta}=\pi^{\alpha\sigma} X_{\sigma\beta}$ or, equivalently,
		$X_{\alpha\beta}=-\pi_{\dt\alpha}^{\bar\tau}A_{\beta\bar\tau}.$
	\end{lemma}
	\begin{proof}
		By definition,
		\[
		sp(n) = \{A \in End(V): g(Ax, y) + g(x, Ay) = 0, AJ_s = J_sA \text{ for } s = 1,2,3\}.
		\]
		
		Setting $g(Ae_\alpha, e_{\bar\beta}) = A_{\alpha \bar\beta}$ for $A \in sp(n)$ we obtain (2).
		To obtain the equivalence between (1) and (2), take an element $A_{\alpha\bar\beta}$ of \eqref{Spn_A} and define $X_{\alpha\beta}=-\pi_{\ \alpha}^{\bar\tau}A_{\beta\bar\tau}.$ Then, 
		\begin{equation*}
			\begin{aligned}
				X_{\beta\alpha}\ =\ -\pi_{\ \beta}^{\bar\tau}A_{\alpha\bar\tau}\ =\ -\pi_{\ \beta}^{\bar\tau}(\mathfrak jA)_{\alpha\bar\tau}\ &=\ -\pi_{\ \beta}^{\bar\tau}\pi^{\bar\sigma}_{\ \alpha}\pi_{\ \bar\tau}^{\mu}A_{\bar\sigma\mu}\\
				&=\ \pi^{\bar\sigma}_{\ \alpha}A_{\bar\sigma\beta}\ =\ -\pi^{\bar\sigma}_{\ \alpha}A_{\beta\bar\sigma}\ =\ X_{\alpha\beta}
			\end{aligned}
		\end{equation*}
		\noindent and thus $X_{\alpha\beta}$ belongs to \eqref{Spn_X}. The converse is similar.
	\end{proof}
	
	The positive definite inner product $g$ on $V$ (cf. \eqref{g-pi}) induces a positive definite  inner product $\lc,\rc$ on $End(V)$ by the formula
	\begin{equation*}
		\lc A,B\rc=\frac{1}{2}\sum_{a=1}^{4n}g(Ah_a,Bh_a),\qquad A,B\in End(V),
	\end{equation*}
	where $\{h_a,\ a=1,\dots 4n\}$ is any $g$-orthonormal basis of $V$. If $\{e_\alpha\}$ is any basis of $W$ (with the conventions from Section~\ref{index_conventions}) and if regarding $sp(n)$ as a subspace of  $End(V)$ (cf. \eqref{Spn_A}), then 
	\begin{equation}\label{killing_form_A}
		\big\lc A,B\big\rc =A^{\alpha\bar\beta}B_{\alpha\bar\beta},\qquad A,B\in sp(n)\subset End(V),
	\end{equation}
	where  $A_{\alpha\bar\beta}=g(Ae_\alpha,e_{\bar\beta})$ and $B_{\alpha\bar\beta}=g(Be_\alpha,e_{\bar\beta})$.
	Alternatively, If regarding $sp(n)$ as a subspace of $\text{S}^2\,{W^*},$  as in \eqref{Spn_X}, then 
	\begin{equation}\label{killing_form_X}
		\big\lc X,Y\big\rc =\pi^{\alpha\gamma}\pi^{\beta\delta}X_{\alpha\beta}Y_{\gamma\delta},\qquad X,Y\in sp(n)\subset S^2W^\ast.
	\end{equation}

	It is easily verified $\big($by using for example formula~(16.3) from \cite{FH}$\big)$ that the so defined inner  product on $sp(n)$ is in fact a scalar multiple of the Killing form.  Explicitly, we have that
	\begin{equation}\label{killing_form_multiple}
		\lc u,v\rc=-\frac{1}{2(n+1)}\text{trace}\Big(\text{ad} (u)\circ \text{ad} (v)\Big),
	\end{equation}
	where $\text{ad}(.)$ is the adjoint representation of $sp(n)$ \big(i.e., $\text{ad}(u)v=[u,v]$, $u,v\in sp(n)$\big).

	A 4-tensor $K$ on $V$ will be called a hyper-K\"ahler curvature-type tensor, if $K\in \Lambda^2\,V^\ast\otimes\Lambda^2\,V^\ast$ and 
	\begin{equation}\label{HK_curv_prop}
		\begin{aligned}
			&K(x,y,z,w)+K(y,z,x,w)+K(z,x,y,w)\ =\ 0,\\
			&K(x,y,J_sz,J_sw)\ =\ K(x,y,z,w), 
		\end{aligned}
	\end{equation}
	for all  $x,y,z,w\in V$ and $s=1,2,3.$
	Notice that the Riemannian curvature tensor of a hyper-K\"ahler manifold \big(i.e., a Riemannian manifold with holonomy $Sp(n)$\big) satisfies these conditions at each point of the manifold \big(cf. \cite{Bes}\big). We shall denote the space of all hyper-K\"ahler curvature-type tensor on $V$ by $\mathfrak R^{Sp(n)}$. 
	\begin{rmrk}\label{rmrk_spnsp1} Observe that it is irrelevant precisely which positively oriented orthonormal basis $J_1,J_2,J_3$ of $\mathcal Q$ $\big($cf. \eqref{spanJ}$\big)$ is used in \eqref{HK_curv_prop}; that is, the subspace $\mathfrak R^{Sp(n)}\subset  (V^\ast)^{\otimes 4}$ is preserved not only under the action of the compact symplectic group $Sp(n)$, but also under the action of the larger group $Sp(n)Sp(1)$ $\big($cf. \eqref{def_spnsp1}$\big)$.
	\end{rmrk}
	Let us denote by $\mathfrak S^{Sp(n)}$ the space of all totally symmetric 4-tensors $$S\ \in\  S^4(W^\ast)+S^4(\overline{W}^*)\ \subset\ S^4(V^{\ast })$$ with the property
	\begin{equation*}
		\overline{S(x,y,z,w)}\ =\ S(\bar x,\bar y,\bar z,\bar w)\ =\ S(J_2x,J_2y,J_2z,J_2w)
	\end{equation*}
	for all $x,y,z,w\in V^{\mathbb C}=W\oplus\overline{W}$. If $\{e_\alpha\}$ is a fixed basis of $W,$ each $S\in \mathfrak S^{Sp(n)}$ produces a totally symmetric array $S_{\alpha\beta\gamma\delta}=S(e_\alpha,e_\beta,e_\gamma,e_\delta)$ (of complex numbers) with the property
	\begin{equation}\label{S_HK}
		S_{\alpha\beta\gamma\delta}=(\mathfrak j S)_{\alpha\beta\gamma\delta}\qquad \Big(cf.\ \text{\eqref{map_j}}\Big).
	\end{equation}
	The converse is also true: each totally symmetric array $S_{\alpha\beta\gamma\delta}$, that satisfies \eqref{S_HK}, represents a unique element  $S\in \mathfrak S^{Sp(n)}$. Clearly, the real vector space $\mathfrak S^{Sp(n)}$ is preserved under the action of $Sp(n)$, but not under the action of $Sp(n)Sp(1)$ (since for its definition we have used the complex structures $J_1$ and $J_2$). 
	
	It is well known  that $\mathfrak R^{Sp(n)}$ and $\mathfrak S^{Sp(n)}$ are isomorphic as representations of $Sp(n)$ (see for example \cite{Sal2}, Proposition~9.3). In  the following lemma we present the respective isomorphism in an explicit form (the proof is straightforward).
	\begin{lemma}\label{1-1 HKC} There is a $Sp(n)$-equivariant isomorphism $\mathcal K:\mathfrak S^{Sp(n)}\rightarrow \mathfrak R^{Sp(n)}$, defined by the formula
		\begin{equation*}
			\mathcal K(S)(x,y,z,w)=S(x,J_2y,z,J_2w)-S(x,J_2y,J_2z,w)
		\end{equation*}
		for all $S\in\mathfrak S^{Sp(n)}$ and $x,y,z,w\in V.$ 
		
		Moreover, the inverse map $\mathcal K^{-1}$ is given by
		\begin{equation*}
			\mathcal K^{-1}(K)(x,y,z,w)=\frac{1}{2}\big(K(x,J_2y,z,J_2w)-K(x,J_3y,z,J_3w)\big)
		\end{equation*}
		for any $K\in\mathfrak R^{Sp(n)}.$ 
		
		Furthermore,
		if $\{e_\alpha\}$ is any basis of $W$ and $K=\mathcal K(S)$, we have
		\begin{equation}\label{def_S_abgd}
			K_{\alpha\bar\beta\gamma\bar\delta}{\ =\ }S_{\alpha\sigma\gamma\tau}\, \pi^{\sigma}_{\ \bar\beta}\pi^{\tau}_{\ \bar\delta},
		\end{equation} 
		where (following the conventions from Section~\ref{index_conventions})
		$K_{\alpha\bar\beta\gamma\bar\delta}=K(e_\alpha,e_{\bar\beta},e_\gamma,e_{\bar\delta})$  and $S_{\alpha\beta\gamma\delta}=S(e_\alpha,e_\beta,e_\gamma,e_\delta)$.
	\end{lemma}

	For a hyper-K\"ahler curvature-type tensor $K\in \mathfrak R^{Sp(n)}$, we can define a linear transformation $\mathcal T_K$ that operates on the Lie algebra $sp(n)$. This transformation maps an element $A\in sp(n)$ to another element in the same algebra. The transformation is defined by the following equation:
	\begin{equation}\label{def_K_0}
		g(\mathcal T_K(A)x,y)=\frac{1}{2}\sum_{a=1}^{4n}K(x,y,h_a,Ah_a).
	\end{equation}
	In this formula, $\{h_a\}_{a=1}^{4n}$ is an arbitrary orthonormal basis of the vector space $V$. A key property of $\mathcal T_K$ is that it is always both traceless and symmetric.

	If $\{e_\alpha\}$ is any basis of $W$,
	\begin{equation}\label{def_K_2}
		A=\mathcal T_K (B)\ \iff\ A^\alpha_{\dt\beta} =K^{\alpha\ \ \ \delta}_{\dt\beta\gamma\dt}\, B^{\gamma}_{\dt\delta}\ \iff\ A_{\alpha\bar\beta} =K_{\alpha\bar\beta\gamma\bar\delta}\, B^{{\gamma\bar\delta}}.
	\end{equation}
	
	If regarding $sp(n)$ as a subspace of $\text{S}^2\,{W^*}$, as in \eqref{Spn_X}, then for any $X=\{X_{\alpha\beta}\}$ and $Y=\{Y_{\alpha\beta}\}$ in $sp(n)$, we have
	\begin{equation}\label{def_K_1}	
		X=\mathcal T_K (Y)\ \iff\ X_{\alpha\beta} \ =\  K_{\alpha\bar\sigma\beta\bar\tau}\, Y^{{\bar\sigma\bar\tau}}.
	\end{equation}

	To show the equivalence between \eqref{def_K_1} and \eqref{def_K_2}, we let $S=\mathcal K^{-1}(K)$ (cf. Lemma~\ref{1-1 HKC})  and  take $X_{\alpha\beta}=-\pi_{\alpha\sigma}A^{\sigma}_{\dt\beta}$, $Y_{\alpha\beta}=-\pi_{\alpha\sigma}B^{\sigma}_{\dt\beta}$. Then, by \ref{def_K_2}, $A=\mathcal T_K (B)$ is equivalent to
	\begin{multline*}
		X_{\alpha\beta}=-\pi_{\alpha\sigma}A^{\sigma}_{\dt\beta}\overset{\eqref{def_K_2}}{\ =\ }
		-\pi_{\alpha\sigma}\hat K^{\sigma\ \ \ \delta}_{\dt\beta\gamma\dt}\, B^{\gamma}_{\dt\delta}
		\ =\ \pi^{\bar\sigma}_{\dt\alpha}\,\pi^{\gamma\mu}\,\pi_{\mu\nu}\,g^{\bar\tau\delta}\,K_{\bar\sigma\beta\gamma\bar\tau}\,B^{\nu}_{\dt\delta}\\
		=\ \pi^{\bar\sigma}_{\dt\alpha}\,\pi^{\gamma\mu}\,g^{\bar\tau\delta}\,K_{\bar\sigma\beta\gamma\bar\tau}\,Y_{\mu\delta}\overset{\eqref{def_S_abgd}}{\ =\ }
		-\pi^{\bar\sigma}_{\dt\alpha}\,\pi^{\gamma\mu}\,g^{\bar\tau\delta}\,\pi^{\xi}_{\dt\bar\sigma}\,\pi^{\zeta}_{\dt\bar\tau}\,S_{\beta\xi\gamma\zeta}\,Y_{\mu\delta}\\
		=\ \pi^{\gamma\mu}\, \pi^{\zeta\delta}\,S_{\beta\alpha\gamma\zeta}\,Y_{\mu\delta}\ 
		= \ K_{\alpha\bar\sigma\beta\bar\tau}\, Y^{{\bar\sigma\bar\tau}}.
	\end{multline*}

	Next, we consider the question: Under what circumstances can an arbitrary linear transformation $L$ on $sp(n)$ be expressed as  $L=\mathcal T_K$ for some $K\in\mathfrak R^{Sp(n)}$? In order to answer this, we shall first construct a certain invariant operator~$\dag$ on $\text{End}\big(sp(n)\big)$ and then show that its 2-eigenspace is precisely what we need.  Let us pick a basis $E_1,\dots,E_m$ of $sp(n)$ \big($m=\dim(sp(n))$\big) and let $E^*_1,\dots E^*_m\in sp(n)$ be the corresponding dual basis; that is, for any $1\le s,t\le m$,
	\begin{equation}\label{def_dual_basis}
		\lc E_s,E_t^*\rc=
		\begin{cases}
			1,\   \text{if}\ s=t\\
			0,\  \text{otherwise}.
		\end{cases}
	\end{equation}
	Using the Lie algebra brackets $[,]$ of $sp(n)$, we define a linear operator $\dag$ on $\text{End}\big(sp(n)\big)$,
	\begin{equation}\label{def_dag}
		(\dag L)X=\sum_{s=1}^m \big[E^*_s,L[E_s,X]\big]\qquad \forall X\in sp(n),\quad \forall L\in \text{End}\big(sp(n)\big).
	\end{equation}
	Clearly, $\dag$ is independent of the choice of the initial basis  $E_1,\dots,E_m$ and by \eqref{killing_form_multiple}, we have
	\begin{equation}\label{dag_id}
		\dag\Big(\text{Id}_{sp(n)}\Big)=-2(n+1)\text{Id}_{sp(n)}.
	\end{equation}
	
	\begin{lemma}\label{T_K} For an element $L\in \text{End}\big(sp(n)\big)$, there exists $K\in\mathfrak R^{Sp(n)}$ so that $L=\mathcal T_K$ iff 
		\begin{equation}\label{dag_eigen}
			\dag L=2L.
		\end{equation}
	\end{lemma}
	\begin{proof}
		By \eqref{Spn_X}, we may regard the complex vector space $\text{S}^2\,{W^*}$ as a complexification of its real subspace $sp(n)$. The Lie bracket in $sp(n)$ induces a bracket operation $[,]$ on  $\text{S}^2\,{W^*}$, 
		\begin{equation*}
			[X,Y]\ =\ \Big\{\pi^{\sigma\tau}\big(X_{\alpha\sigma}Y_{\tau\beta}+X_{\beta\sigma}Y_{\tau\alpha}\big)\Big\}
		\end{equation*}
		for any  $X=\{X_{\alpha\beta}\}$ and $Y=\{Y_{\alpha\beta}\}$ in $\text{S}^2\,{W^*}$.
		
		Let us pick a basis $\sharp^{\alpha\beta}\overset{def}{\ \ = \ \ }\frac{1}{2}(e^\alpha\otimes e^\beta+e^\beta\otimes e^\alpha)$, $1\le \alpha\le \beta\le 2n$ for $\text{S}^2\,{W^*}$, where $\{e^\alpha\}$ is a basis of $W^\ast$.  Using the invariant inner product \eqref{killing_form_X}, we calculate (cf. \eqref{def_dual_basis}) that the corresponding dual basis 
		$\sharp^\ast_{\alpha\beta}$ ($1\le \alpha\le \beta\le 2n$) is given by 
		$$\sharp^\ast_{\alpha\beta}=
		\begin{cases}
			\$_{\alpha\alpha},\   \text{if}\ \alpha=\beta\\
			2\$_{\alpha\beta},\  \text{otherwise},
		\end{cases}$$
		where $\$_{\alpha\beta}\overset{def}{\ =\ }\sum_{\sigma,\tau=1}^{2n}\pi_{\alpha\sigma}\pi_{\beta\tau}\,\sharp^{\sigma\tau}$. 
		A straightforward calculation shows that 
		\begin{equation}\label{Lie_br_sp2}
			\begin{array}{ll}
				\big[\sharp^{\alpha\beta},\sharp^{\gamma\delta}\big]=\frac{1}{2}\big(\pi^{\alpha\gamma}\sharp^{\beta\delta}+\pi^{\beta\delta}\sharp^{\alpha\gamma}+\pi^{\alpha\delta}\sharp^{\beta\gamma}+\pi^{\beta\gamma}\sharp^{\alpha\delta}\big),\\[7pt]
				
				\big[\$_{\alpha\beta},\$_{\gamma\delta}\big]=\frac{1}{2}\big(\pi_{\alpha\gamma}\$_{\beta\delta}+\pi_{\beta\delta}\$_{\alpha\gamma}+\pi_{\alpha\delta}\$_{\beta\gamma}+\pi_{\beta\gamma}\$_{\alpha\delta}\big).
			\end{array}
		\end{equation}

		A linear transformation $L$ on the Lie algebra $sp(n)$ can be represented by a unique set of components $S_{\alpha\beta\gamma\delta}$, in a given basis. This array of complex numbers possesses a key symmetry: it is unchanged if you swap the indices $\alpha$ and $\beta$, or if you swap the indices $\gamma$ and $\delta$. We have also that 
		\mbox{$
			L(\$_{\alpha\beta})=\sum_{\gamma,\delta=1}^{2n}S_{\gamma\delta\alpha\beta}\,\sharp^{\gamma\delta}
			$}. Using this, we calculate 
		\begin{equation}
			\begin{aligned}
				\label{dag_identity}
				(\dag L)(\$_{\alpha\beta})\ &=\ \sum_{\sigma,\tau=1}^{2n}\Big[\sharp^{\sigma\tau},L[\$_{\sigma\tau},\$_{\alpha\beta}]\Big]\\ 
				&=\ \frac{1}{2}
				\Big[\sharp^{\sigma\tau},
				L\big(\pi_{\sigma\alpha}\$_{\tau\beta}+\pi_{\tau\beta}\$_{\sigma\alpha}+\pi_{\sigma\beta}\$_{\tau\alpha}+\pi_{\tau\alpha}\$_{\sigma\beta}\big)\Big]\\
				&=\ \frac{1}{2}\Big(\pi_{\sigma\alpha}S_{\mu\nu\tau\beta}+\pi_{\tau\beta}S_{\mu\nu\sigma\alpha}+\pi_{\sigma\beta}S_{\mu\nu\tau\alpha}+\pi_{\tau\alpha}S_{\mu\nu\sigma\beta}\Big)\big[\sharp^{\sigma\tau},\sharp^{\mu\nu}\big]\\
				&=\ \frac{1}{2}\Big\{S_{\alpha\gamma\delta\beta}+S_{\beta\gamma\delta\alpha}+S_{\alpha\delta\gamma\beta}+S_{\beta\delta\gamma\alpha}\\
				&\hspace{2cm}-\ \pi^{\sigma\tau}\big(\pi_{\gamma\alpha}S_{\delta\sigma\tau\beta}+\pi_{\delta\alpha}S_{\gamma\sigma\tau\beta}+\pi_{\gamma\beta}S_{\delta\sigma\tau\alpha}+\pi_{\delta\beta}S_{\gamma\sigma\tau\alpha}\big)\Big\}\sharp^{\gamma\delta}.
			\end{aligned}
		\end{equation}

		If we assume that $L=\mathcal T_K$, then by Lemma~\ref{1-1 HKC},  the array $S_{\alpha\beta\gamma\delta}$ must be totaly symmetric and  equation \eqref{dag_identity} yields \eqref{dag_eigen}.
		
		Conversely, assume that \eqref{dag_eigen} is satisfied and consider the two traces  $T_{\alpha\beta}$  and $l$,
		\begin{equation}
			T_{\alpha\beta}\overset{def}{\ =\ }\pi^{\sigma\tau}S_{\alpha\sigma\tau\beta},\qquad l\overset{def}{\ =\ }\pi^{\sigma\tau}T_{\sigma\tau}.
		\end{equation}
		We calculate
		\begin{equation*}
			2T_{\delta\beta}\ =\ 2\pi^{\sigma\tau}S_{\delta\gamma\tau\beta}\overset{\eqref{dag_eigen}}{\ =\ }\pi^{\sigma\tau}(\dag S)_{\delta\gamma\tau\beta}
			\overset{\eqref{dag_identity}}{\ =\ }\frac{1}{2}\Big\{T_{\beta\delta}-(2n+3) T_{\delta\beta}-l\pi_{\delta\beta}\Big\}
		\end{equation*}
		and thus $T_{\alpha\beta}=0$, $l=0$. By \eqref{dag_identity}, 
		\begin{equation}\label{dag_eig_sec}
			2S_{\alpha\beta\gamma\delta}=S_{\alpha\gamma\delta\beta}+S_{\beta\gamma\delta\alpha}+S_{\alpha\delta\gamma\beta}+S_{\beta\delta\gamma\alpha}.
		\end{equation}
		The symmetries on the right side of equation \eqref{dag_eig_sec} demonstrate that $S_{\alpha\beta\gamma\delta}=S_{\gamma\delta\alpha\beta}$. This means that equation  \eqref{dag_eig_sec} can be restated as
		$
		3S_{\alpha\beta\gamma\delta}=S_{\alpha\beta\gamma\delta}+S_{\beta\gamma\alpha\delta}+S_{\gamma\alpha\beta\delta}.
		$
		The right side of this new equation is clearly symmetric with respect to the indices $\alpha,\beta,\gamma$. Therefore, the array $S_{\alpha\beta\gamma\delta}$ must be totally symmetric. On the other hand, by assumption, $L$ is an endomorphism of $\text{S}^2\,{W^*}$ that preserves its real structure, since it preserves $sp(n)$ (which we identify with the real part of $\text{S}^2\,{W^*}$), i.e., we have 
		$$(\mathfrak j S)_{\alpha\beta\gamma\delta}=S_{\alpha\beta\gamma\delta}.$$ 
		Therefore, $K=\mathcal K(S)$ is a hyper-K\"ahler curvature-type tensor and $L=\mathcal T_K$.
	\end{proof}

	\section{A geometric interpretation for the discriminant of a cubic polynomial}
	This section demonstrates that Formula \eqref{cubic_dis} for the discriminant of a cubic polynomial determines a specific $Sp(2)Sp(1)$-orbit, $\mathcal C$, within the space of all hyper-K\"ahler curvature-type tensors (cf. \eqref{HK_curv_prop}), $\mathfrak R^{Sp(2)}$, associated to any 8-dimensional vector space equipped with a hyper-Hermitian structure. We shall provide a basis-independent description of $\mathcal C$ and prove that it is isomorphic to the homogeneous space
	$$\frac{Sp(2)Sp(1)}{G\,Sp(1)},$$ where $G$ is the image of $Sp(1)\cong SU(2)$, determined by its unique (up to an isomorphism) 4-dimensional complex (or 8-dimensional real) irreducible representation. The product group $G\,Sp(1)$ is isomorphic to $SO(4)=\frac{Sp(1)\times Sp(1)}{\mathbb Z_2}$ and corresponds to an irreducible representation of $SO(4)$ on $\mathbb R^8$.
	
	Explicitly, the group $G$---for which we shall use the notation $G=Sp(1)_{\text{ir}}$---is constructed as follows:
	Consider the standard representation of $Sp(1)$ on $\Delta=\mathbb C^2$ and let $W$ be the space of all symmetric rank-$3$ (contravariant) tensors on $\Delta$, 
	\begin{equation}\label{W-sym}
		W=\text{S}^3\Delta.
	\end{equation}  We may think of $W$ as the space of polynomials in one variable of degree at most three (with complex coefficients) and therefore, the arising 4-dimensional (complex) representation of $Sp(1)$ on $W$ is irreducible. By definition, $Sp(1)_{\text{ir}}$ is the image of $Sp(1)$ in the general linear group $GL_{\mathbb C}(W)$ determined by this representation. There is a natural hyper-Hermitian structure on $W$ with respect to which we have an inclusion 
	\begin{equation}\label{incl_ir_sp2}
		Sp(1)_{\text{ir}}\subset Sp(2).
	\end{equation}
	Indeed, the complex vector space $\Delta$ possesses a natural 2-form $\pi\in \Lambda^2\Delta^\ast$ and an anti-linear transformation $j$ satisfying \eqref{q-str}, that are $Sp(1)$-equvariant. The tensor products $\pi\otimes\pi\otimes\pi$ and $j\otimes j\otimes j$ determine a 2-form and a compatible anti-linear transformation on $W$---to be denoted also by $\pi$ and $j$ (with a slight abuse of notation)---that satisfy \eqref{q-str} and thus define an $Sp(1)_{\text{ir}}$-equivariant hyper-Hermitian structure on $W$. The stabilizer of the pair $(\pi,j)$ in the general linear group $GL_{\mathbb C}(W)$ is $Sp(2)$  and we have the inclusion \eqref{incl_ir_sp2}.

	If considering the underlying 8-dimensional real vector space $V,$ then $V\otimes_{\mathbb R} \mathbb C=W\oplus \overline{W}$ and we are in a situation similar to this from the discussion in Section~\ref{preliminaries} (cf. ~\eqref{splitVc} and~\eqref{g-pi}). Using the same construction and notation, we have an induced positive definite inner product $g$ on $V$ and a compatible hyper-Hermitian structure $J_1,J_2,J_3$ with respect to which
	$$Sp(1)_{\text{ir}}\subset Sp(2)\subset GL(V).$$

	\subsection{Representation in coordinates}
	Let us fix a basis $d_1, d_2=j(d_1)$ of $\Delta$, so that $\pi(d_1,d_2)=1$.  The abstract Lie algebra $sp(1)$ has a basis $\hat{\mathcal  E}_1, \hat{\mathcal E}_2,\hat{\mathcal E}_3$, satisfying
	$[\hat{\mathcal E}_1,\hat{\mathcal E}_2]=\hat{\mathcal E}_3$, $[\hat{\mathcal E}_2,\hat{\mathcal E}_3]=\hat{\mathcal E}_1$ and $[\hat{\mathcal E}_3,\hat{\mathcal E}_1]=\hat{\mathcal E}_2$, that acts on $\Delta$ via the complex matrices
	\begin{equation}\label{rep_delta}
		\hat{\mathcal E}_1\mapsto\frac{1}{2}\begin{pmatrix}
			-i &\ \, 0\\
			\ \,0 & \ \,i
		\end{pmatrix},
		\qquad\hat{\mathcal E}_2\mapsto
		\frac{1}{2}\begin{pmatrix}
			\ \, 0 & -1\\
			\ \,1 &\ \, 0
		\end{pmatrix},
		\qquad\hat{\mathcal E}_3\mapsto
		\frac{1}{2}\begin{pmatrix}
			0 & i\\
			i & 0
		\end{pmatrix}.
	\end{equation}

	Since $W=\text{S}^3\Delta$, we can take an $Sp(2)$-adapted basis $\hat e_1,\,\hat e_2,\,\hat e_3,\,\hat e_4$ for $W$ (cf. \eqref{constants}) by setting:
	\begin{equation}\label{basis_hat}
		\begin{split}
			\hat e_1&=d_1^{\otimes 3},\quad \hat e_2=\frac{1}{\sqrt{3}}\Big(d_1\otimes d_1\otimes d_2+d_1\otimes d_2\otimes d_1+d_2\otimes d_1\otimes d_1\Big),\\
			\hat e_3&=d_2^{\otimes 3},\quad \hat e_4=-\frac{1}{\sqrt{3}}\Big(d_1\otimes d_2\otimes d_2+d_2\otimes d_1\otimes d_2+d_2\otimes d_2\otimes d_1\Big).
		\end{split}
	\end{equation} 
	With respect to this basis the action of $sp(1)_{\text{ir}}$ on $W$ is given by the mapping $\hat{\mathcal E}_s\mapsto E_s$, $s=1,2,3$, where
	
	\begin{equation}\label{rep_W}
		E_1=\text{\scalebox{0.87}{$\left(
				\begin{smallmatrix}
					-\frac{3i}{2} &\ \, 0&\ \,0 &\ \,0\ \\
					\ \,0 & -\frac{i}{2}&\ \,0 &\ \,0\ \\
					\ \,0 &\ \, 0&\ \,\frac{3i}{2} &\ \,0\ \\
					\ \,0 &\ \, 0&\ \,0 &\ \,\frac{i}{2}\ \\
				\end{smallmatrix}\right)$}}\quad
		E_2=\text{\scalebox{0.87}{$\left(
				\begin{smallmatrix}
					\ \,0 &-\frac{\sqrt{3}}{2}&\ \,0 &\ \,0\\
					\ \,\frac{\sqrt{3}}{2} &\ \, 0&\ \,0 &\ \,1\\
					\ \,0 &\ \, 0&\ \,0 &-\frac{\sqrt{3}}{2}\\
					\ \,0 &-1&\ \,\frac{\sqrt{3}}{2} &\ \,0\\
				\end{smallmatrix}\right)$}}\quad
		E_3= \text{\scalebox{0.87}{$\left(
				\begin{smallmatrix}
					\ \,0 &\ \, \frac{i\sqrt{3}}{2}&\ \,0 &\ \,0\\
					\frac{i\sqrt{3}}{2} &\ \, 0&\ \,0 & -i\\
					\ \,0 &\ \, 0&\ \,0 &-\frac{i\sqrt{3}}{2}\\
					\ \,0 &-i& -\frac{i\sqrt{3}}{2} &\ \,0\\
				\end{smallmatrix}
				\right)$}}
	\end{equation}
	
	\vspace{0.3cm}
	
	Consider the orthogonal projection
	\begin{equation}\label{pr_def}
		\hat{\mathcal P}:sp(2)\rightarrow sp(1)_{\text{ir}}.
	\end{equation}
	A straightforward calculation shows that (cf. \eqref{def_dag}) 
	\begin{equation*}
		\dag \hat{\mathcal P}=2\hat{\mathcal P}-\frac{12}{5}\text{Id}.
	\end{equation*}
	Since $\dag \text{Id}=-6\text{Id}$ (cf. \eqref{dag_id}), we have that $\dag \big(\hat{\mathcal P}-\frac{3}{10}\text{Id}\big)=2\big(\hat{\mathcal P}-\frac{3}{10}\text{Id}\big)$ and therefore, by Lemma~\ref{T_K}, there exists a $\hat K\in\mathfrak R^{Sp(2)},$  so that
	\begin{equation}\label{proj_prop}
		5\Big(\hat{\mathcal P}-\frac{3}{10}\text{Id}\Big)=\mathcal T_{\hat K}.
	\end{equation}
	
	Our main task here is to calculate explicitly the components $\hat S_{\alpha\beta\gamma\delta}$ of the totally symmetric 4-tensor $\hat S=\mathcal K^{-1}(\hat K)$ (cf. Lemma~\ref{1-1 HKC}) with respect to the basis \eqref{basis_hat}. Before doing this, however, we shall collect some further technical data about the induced structure on $W$. It will  be useful to consider the three symmetric matrices $\Upsilon_1$, $\Upsilon_2$ and $\Upsilon_3$ obtained by multiplying the matrices \eqref{rep_W} from the left with
	\begin{equation*}\left(
		\begin{smallmatrix}
			\, 0 &\ \, 0&-1 &\ \,0\\
			\,0 &\ \, 0&\ \,0 &-1\\
			\,1 &\ \, 0&\ \,0 &\ \,0\\
			\,0 &\ \, 1&\ \,0 &\ \,0\\
		\end{smallmatrix}\right),
	\end{equation*}
	i.e.,
	\begin{equation}\label{upsilon}
		\Upsilon_1= \text{\scalebox{0.87}{$\left(
				\begin{smallmatrix}
					\ \, 0 &\ \, 0&-\frac{3i}{2} &\ \,0\ \\
					\ \,0 &\ \,0&\ \,0 &-\frac{i}{2}\ \\
					-\frac{3i}{2} &\ \, 0&\ \, 0 &\ \,0\ \\
					\ \,0 &-\frac{i}{2}&\ \,0 &\ \,0\ \\
				\end{smallmatrix}\right)$}}\quad
		\Upsilon_2= \text{\scalebox{0.87}{$\left(
				\begin{smallmatrix}
					\ \,0 &\ \, 0&\ \,0 &\ \,\frac{\sqrt{3}}{2}\\
					\ \,0 &\ \, 1&-\frac{\sqrt{3}}{2} &\ \,0\\
					\ \,0 &-\frac{\sqrt{3}}{2}&\ \,0 &\ \, 0\\
					\ \,\frac{\sqrt{3}}{2} &\ \, 0&\ \,0 &\ \,1\\
				\end{smallmatrix}\right)$}}\quad
		\Upsilon_3= \text{\scalebox{0.87}{$\left(
				\begin{smallmatrix}
					\ \,0 &\ \, 0&\ \,0 &\ \,\frac{i\sqrt{3}}{2}\\
					\ \,0 &\ \, i&\ \,\frac{i\sqrt{3}}{2} &\ \,0\\
					\ \,0 &\ \, \frac{i\sqrt{3}}{2}&\ \,0 &\ \, 0\\
					\ \,\frac{i\sqrt{3}}{2} &\ \, 0&\ \,0 &-i\\
				\end{smallmatrix}
				\right)$}}.
	\end{equation}
	We shall denote by $(\Upsilon_s)_{\alpha\beta}$ the entry in the $\alpha$-row and $\beta$-column of $\Upsilon_s$. Then, the action \eqref{rep_W} of $\hat{\mathcal E}_s$ on $W$ can be written in the form 
	\begin{equation}\label{upsilon_E}
		\hat{\mathcal E}_s\mapsto\pi^{\alpha\sigma}(\Upsilon_s)_{\sigma\beta},\qquad s=1,2,3.
	\end{equation}
	\begin{lemma}\label{upsilon_lemma}
		We have the following identities:
		
		\vspace{0.3cm}\par (1)\hspace{0.3cm}$(\Upsilon_s)_{\alpha\beta}=(\Upsilon_s)_{\beta\alpha}$, $(\mathfrak j\Upsilon_s)_{\alpha\beta}\ =\ (\Upsilon_s)_{\alpha\beta}$;
		
		\vspace{0.3cm}\par (2)\hspace{0.3cm}$\pi^{\alpha\sigma}\,\pi^{\beta\tau}\,(\Upsilon_s)_{\alpha\beta}\,(\Upsilon_t)_{\sigma\tau}\ =\ 
		\begin{cases}
			5,\ \text{if}\ s=t,\\
			0,\ \text{otherwise}.
		\end{cases}$;
		
		\vspace{0.3cm}\par (3)\hspace{0.3cm}$\pi^{\sigma\tau}\Big((\Upsilon_i)_{\alpha\sigma}(\Upsilon_j)_{\beta\tau}+(\Upsilon_i)_{\beta\sigma}(\Upsilon_j)_{\alpha\tau}\Big)\ =\ (\Upsilon_k)_{\alpha\beta}$;\ \footnotemark \footnotetext{$(ijk)$ is any positive permutation of $(123)$}
		
		\vspace{0.3cm}\par (4)\hspace{0.3cm}$\pi^{\sigma\tau}\,\Big((\Upsilon_s)_{\alpha\sigma}\,\hat S_{\tau\beta\gamma\delta}+(\Upsilon_s)_{\beta\sigma}\,\hat S_{\tau\alpha\gamma\delta}+(\Upsilon_s)_{\gamma\sigma}\,\hat S_{\tau\alpha\beta\delta}+(\Upsilon_s)_{\delta\sigma}\,\hat S_{\tau\alpha\beta\gamma}\Big)\ =\ 0$;
		
		\vspace{0.3cm}\par (5)\hspace{0.3cm}$\sum_s(\Upsilon_s)_{\alpha\beta}(\Upsilon_s)_{\gamma\delta}\ =\ \hat S_{\alpha\beta\gamma\delta}+\frac{3}{4}\big(\pi_{\alpha\gamma}\pi_{\beta\delta}+\pi_{\alpha\delta}\pi_{\beta\gamma}\big)$.
		
		\vspace{0.3cm}\par (6)\hspace{0.3cm}$(\Upsilon_s)_{\sigma\tau}\,
		\pi^{\sigma\gamma}\,\pi^{\tau\delta}\,\hat S_{\alpha\beta\gamma\delta} 
		\ =\ \frac{7}{2}\,(\Upsilon_s)_{\alpha\beta},\quad \forall s=1,2,3$.
		
		\vspace{0.3cm}\par (7)\hspace{0.3cm}$\pi^{\sigma\tau}\,\sum_s(\Upsilon_s)_{\alpha\sigma}\,(\Upsilon_s)_{\beta\tau} \ =\ \frac{15}{4}\,\pi_{\alpha\beta}$.
	\end{lemma}
	\begin{proof}
		Since, by construction, each $E_s$ is an element of $sp(2), $ {\it(1)} follows from \eqref{Spn_X}. The relation {\it(2)} is obtained by calculating the trace for each of the product matrices $E_s E_t$.  {\it(3)}~is equivalent to the relation  $\big[ E_i,  E_j\big]=  E_k$. Since each $E_s$ belongs to the Lie algebra of the stabilizer group of $\hat S$, we must have 
		\begin{equation*}
			\hat S(  E_s\cdot x,y,z,w)+\hat  S(x, E_s\cdot y,z,w)+\hat  S( x,y, E_s\cdot z,w)+ \hat  S(x,y,z, E_s\cdot w)=0,
		\end{equation*} 
		for any $x,y,z,w \in W$, and hence {\it(4)} is satisfied. {\it(5)} is equivalent to \eqref{proj_prop}.
		To show{\it(6)} we calculate using {\it(2)} and {\it(5)} that
		\begin{equation*}
			\begin{aligned}
				(\Upsilon_s)_{\sigma\tau}\,
				\pi^{\sigma\gamma}\,\pi^{\tau\delta}\,\hat S_{\alpha\beta\gamma\delta} 
				\ &=\ (\Upsilon_s)_{\sigma\tau}\,
				\pi^{\sigma\gamma}\,\pi^{\tau\delta}\,\Big(\sum_t(\Upsilon_t)_{\alpha\beta}(\Upsilon_t)_{\gamma\delta}\ -\ \frac{3}{4}\big(\pi_{\alpha\gamma}\pi_{\beta\delta}+\pi_{\alpha\delta}\pi_{\beta\gamma}\big)\Big)\\
				&=\ 5\,(\Upsilon_s)_{\alpha\beta}\ -\ \frac{3}{4}\,(\Upsilon_s)_{\alpha\beta}\ -\ \frac{3}{4}\,(\Upsilon_s)_{\beta\alpha}\ =\  \frac{7}{2}\,(\Upsilon_s)_{\alpha\beta}.
			\end{aligned}
		\end{equation*}
		{\it(7)} is a direct consequence of {\it(5)}.
	\end{proof}

	Based on classical invariant theory, it is a well-established fact that the subalgebra of $Sp(1)$-invariants within $S^4W^*$ is one-dimensional.
	This means there is only one fundamental, non-trivial invariant of this type. Furthermore, it is a known result that this single invariant can be represented by the discriminant of a cubic polynomial. Therefore, after an appropriate change of basis for the vector space $W,$ the 4-tensor $\hat S$ will be equivalent to this discriminant.
	We now turn to the task of computing the components $\hat S_{\alpha\beta\gamma\delta}$ of $\hat S$ with respect to the basis \eqref{basis_hat}. Using equation {\it(5)} from Lemma~\ref{upsilon_lemma}, we calculate
	\begin{equation}\label{new-dis}
		\begin{aligned}
			&\hat{S}_{\alpha\beta\gamma\delta}\,x^\alpha x^\beta x^\gamma x^\delta\  =\ \sum_s\Big(\sum_{\alpha,\beta}\Upsilon^s_{\alpha\beta}\,x^\alpha\,x^\beta\Big)^\pow{2}\\
			&\hspace{0.5cm}=\ \Big(-3\,i\,x^\ind{1}\,x^\ind{3}-i\,x^\ind{2}\,x^\ind{4}\Big)^\pow{2}\ 
			+\ \Big(\sqrt{3}\,x^\ind{1}\,x^\ind{4}+(x^\ind{2})^\pow{2}-\sqrt{3}\,x^\ind{2}\,x^\ind{3}+(x^\ind{4})^\pow{2}\Big)^\pow{2}\\
			&\hspace{5cm}+\ \Big(i\,\sqrt{3}\,x^\ind{1}\,x^\ind{4}+i\,(x^\ind{2})^\pow{2}+i\,\sqrt{3}\,x^\ind{2}\,x^\ind{3}-i\,(x^\ind{4})^\pow{2}\Big)^\pow{2}\\[0.3cm]
			&\hspace{0.5cm}=-18\, x^\ind{1} x^\ind{2} x^\ind{3} x^\ind{4}+4\sqrt{3}\,x^\ind{1}(x^\ind{4})^\pow{3}
			-4\sqrt{3}\,(x^\ind{2})^\pow{3}x^\ind{3}+3(x^\ind{2})^\pow{2}(x^\ind{4})^\pow{2}-9(x^\ind{1})^\pow{2}(x^\ind{3})^\pow{2}.
		\end{aligned}
	\end{equation}
	
	Observe that, by the simple substitution 
	$$x^\ind{1}=a,\ x^\ind{2}=\frac{b}{\sqrt{3}},\ x^\ind{3}=d,\ x^\ind{4}=-\frac{c}{\sqrt{3}},$$
	the expression in the last line of \eqref{new-dis} (multiplied by $3$)  produces Formula~\eqref{cubic_dis} for the discriminant of a cubic polynomial. 
	
	\begin{rmrk}\label{rmk-cd} Although the classical Formula \eqref{cubic_dis} for the discriminant looks a bit simpler than the last line in \eqref{new-dis}, we shall prefer using the latter, since it has the advantage of being obtained on an $Sp(2)$-adapted basis of $W$.
	\end{rmrk}

	\subsection{Quaternionic characterization of  the formula for the discriminant} \label{sp2-inv}
	Let us now reverse our perspective by starting with an arbitrary 8-dimensional vector space $V,$ equipped with a hyper-Hermitian structure for which we shall keep the usual notation $g,\, J_1,J_2$ and $J_3.$ We have also an induced representation of $Sp(2)Sp(1)$ on $\mathfrak R^{Sp(2)}$---the space of all hyper-K\"ahler curvature-type tensors on $V$. The $i$-eigenspace $W\subset V^{\mathbb{C}}$ of the complex structure $J_1$ is a 4-dimensional complex subspace and we have $V^{\mathbb C}=W\oplus \overline{W}.$

	\begin{dfn}\label{dfn_cubic_dis} We define $\mathcal C\subset \mathfrak R^{Sp(2)}$ to be the set of all $K\in \mathfrak R^{Sp(2)}$ for which $S=\mathcal K^{-1}(K)$  $($cf. Lemma~\ref{1-1 HKC}$)$ is given by the formula  for the discriminant of a cubic polynomial $($cf. \eqref{new-dis}$)$, i.e., 
		\begin{multline}\label{def-new-dis}
			S(x,x,x,x)\ =\  {S}_{\alpha\beta\gamma\delta}\,x^\alpha x^\beta x^\gamma x^\delta\\
			=-18\, x^\ind{1} x^\ind{2} x^\ind{3} x^\ind{4}+4\sqrt{3}\,x^\ind{1}(x^\ind{4})^\pow{3}
			-4\sqrt{3}\,(x^\ind{2})^\pow{3}x^\ind{3}+3(x^\ind{2})^\pow{2}(x^\ind{4})^\pow{2}-9(x^\ind{1})^\pow{2}(x^\ind{3})^\pow{2}
		\end{multline}
		for all $x=x^\alpha e_\alpha\in W,$ w.r.t. some $Sp(2)$-adapted basis  $\{e_\alpha\}$ of $W$.  
	\end{dfn}

	Clearly, $\mathcal C$ is a distinguished $Sp(2)$-orbit in $\mathfrak R^{Sp(2)}$.
	Our aim here is to find a set of coordinate-free conditions that characterize the elements of $\mathcal C$. We do this by exploring the algebraic properties of the Lie algebra bracket $[,]$ in $sp(2)$.

	\begin{thrm}\label{sp2-characterization} An element $K\in \mathfrak R^{Sp(2)}$ belongs to the orbit $\mathcal C$  iff the endomorphism $\mathcal T_K$ of $sp(2)$  \big(cf. \eqref{def_K_0}\big) satisfies the following two conditions:
		
		\vspace{0.1cm}
		
		(I) $\ $ $(2\mathcal T_K-7\,\text{Id})(2\mathcal T_K+3\,\text{Id})\ =\ 0.$
		
		\vspace{0.1cm}
		
		(II) $\ $ For all $A,B\in sp(2)$,
		$$\big[\mathcal T_KA,\mathcal T_KB\big]-\mathcal T_K\big[\mathcal T_KA, B\big]\ =\ \frac{3}{2}\,\Big(\mathcal T_K\big[A,B\big]-\big[A,\mathcal T_KB\big]\Big).$$

	\end{thrm}

	\begin{proof} For any  $K\in \mathfrak R^{Sp(2)}$, consider
		the endomorphism $\mathcal P$ of $sp(2)$,
		\begin{equation}\label{def_P}
			\mathcal P\overset{def}{\ =\ }\frac{1}{5}\mathcal T_K+\frac{3}{10}\text {Id}.
		\end{equation}
		It is easily seen that (I) is equivalent to 
		\begin{equation}\label{I_p2}
			\mathcal P^2=\mathcal P,
		\end{equation} 
		whereas (II) is equivalent to
		\begin{equation}\label{II_[pp]}
			\Big[\mathcal PA,\mathcal P B\Big]=\mathcal P\Big[\mathcal PA,B\Big].
		\end{equation}

		Consider again the orthogonal projection  $\hat{\mathcal P}$ from~ \eqref{pr_def}.  Observe that if an endomorphism $X\in sp(2)$ is orthogonal to $sp(1)_{\text{ir}}\subset sp(2)$, then so is the commutator $[\mathcal E,X]$ for every  $\mathcal E\in sp(1)_{\text{ir}}$. This means that we have
		\begin{equation}\label{pr_property_1}
			\hat{\mathcal P}\Big[\hat{\mathcal P}(A),B\Big]=\Big[\hat{\mathcal P}(A),\hat{\mathcal P}(B)\Big],\qquad \forall A,B\in sp(2).
		\end{equation}

		If $K\in\mathcal C$ then $S=\mathcal K^{-1}(K)$ is given by \eqref{def-new-dis} with respect to some $Sp(2)$-adapted basis $\{e_\alpha\}$ of $W$. Let $\varphi$ to the linear map that sends each basis vector $\hat e_\alpha$, as defined in \eqref{basis_hat}, to  $e_\alpha$. Then, $\varphi\in Sp(2)$ and it maps $\mathcal T_{\hat K}$ to $\mathcal T_{K}$. Therefore, in view of \eqref{proj_prop}, $\varphi$ maps also the endomorphism $\hat{\mathcal P}$ to the endomorphism $\mathcal P$ defined by \eqref{def_P}. The properties (I) and (II) of $\mathcal P$  follow from $\hat{\mathcal P}^2=\hat{\mathcal P}$ and \eqref{pr_property_1}.
		
		Conversely, let $K\in \mathfrak R^{Sp(2)}$ be any element for which the corresponding endomorphism $\mathcal P$ (defined by \eqref{def_P}) satisfies conditions (I) and (II). Then, by \eqref{I_p2}, $\mathcal P$ is a projection and its image is a certain subspace 
		$Q\subset sp(2).$ The endomorphism $\mathcal T_K$  has eigenvalues $\frac{7}{2}$ and $-\frac{3}{2}$ with corresponding eigenspaces  $Q$ and $Q^{\perp}$. Since $\mathcal T_K$ is a traceless operator, we must have \big(noting that the dimension of $sp(2)$ is~10\big)
		$$0=\text{trace}\ \mathcal T_K=\frac{7}{2} (\dim Q)-\frac{3}{2}(10-\dim Q)$$
		and consequently $\dim Q=3$. By \eqref{II_[pp]}, $Q$ is a Lie subalgebra of $sp(2)$ with an invariant inner product $\lc,\rc$. Each 3-dimensional real Lie algebra with an invariant positive-definite inner product is either abelian or isomorphic to $sp(1).$ 
		Since the maximal possible dimension of an abelian subalgebra of $sp(2)$ is $2$, $Q$ cannot be abelian, hence we must have $Q\cong sp(1)$. All possible embeddings of the Lie algebra $sp(1)$ into $sp(2)$ are well understood; they all come from a respective representation of the Lie group $Sp(1)$ on the complex vector space $W$ and only one of them (up to an isomorphism) is irreducible. In order to show that $Q$  corresponds precisely to this irreducible representation, we assume the opposite:  Suppose that under the action of $Q$, the vector space $W$ splits, $W=W_1\oplus W_2$, with $W_1$ and $W_2$ being 2-dimensional $Q$-invariant subspaces. Notice that since $\dag \mathcal T_K=2\mathcal T_K$ (cf. Lemma~\ref{T_K}) and $\dag\, \text{Id}=-6\,\text{Id}$,  \eqref{def_P} yields 
		\begin{equation*}
			\mathcal P=\frac{1}{2}\dag \mathcal P +\frac{6}{5}\text{Id}
		\end{equation*}
		and 
		\begin{equation}\label{def_P_cont}
			25(\dag\mathcal P)^2+70\,\dag \mathcal P +24\,\text{Id}\ =\ 0.
		\end{equation}

		Let us first consider the case where  the action of $Q$ on both factors $W_1$ and $W_2$ is non-trivial.  We can pick a $Sp(2)$-adapted frame $\{e_\alpha\}$ of $W$ so that $e_1,e_3\in W_1$, $e_3,e_4\in W_2$ and the action of $Q$ on $W$ is given by the representation
		\begin{equation}\label{reduc_nd_case}
			\hat{\mathcal E}_1\mapsto\left(
			\begin{smallmatrix}
				\ \,0 &\ \,0 &\ \, \frac{1}{2} &\ \,0\\
				\ \, 0&\ \, 0&\ \,0 &\ \, \frac{1}{2}\\
				-\frac{1}{2} &\ \, 0&\ \,0 &\ \,0\\
				\ \,0 &-\frac{1}{2} & \ \,0 &\ \,0\\
			\end{smallmatrix}\right)\quad
			\hat{\mathcal E}_2\mapsto\left(
			\begin{smallmatrix}
				\ 0 &\ 0 &  \frac{i}{2} & \ 0\\
				\ 0 &\ 0 &\  0 & \  \frac{i}{2}\\
				\ \frac{i}{2} & \ 0 &\  0  & \ 0\\
				\ 0 &\ \frac{i}{2} &\  0 & \ 0\\
			\end{smallmatrix}\right)\quad
			\hat{\mathcal E}_3\mapsto\left(
			\begin{smallmatrix}
				\ \, \frac{i}{2} &\ \, 0&\ \,0 &\ \,0\\
				\ \,0 &\ \, \frac{i}{2}&\ \,0 &\ \,0\\
				\ \,0 &\ \, 0&\, -\frac{i}{2} &\ \,0\\
				\ \,0 &\ \, 0&\ \,0 &\,-\frac{i}{2}\\
			\end{smallmatrix}\right),\\
		\end{equation}
		where $\hat{\mathcal E}_1,\hat{\mathcal E}_2,\hat{\mathcal E}_3$ is the standard bases of $Q\cong sp(1)$ from \eqref{rep_delta}.
		
		For the orthogonal projection $\mathcal P:sp(2)\rightarrow Q$, in this case, we obtain, after some straightforward calculations,  
		\begin{equation*}
			(\dag \mathcal P)^2=-2\,\dag \mathcal P,
		\end{equation*}
		which clearly contradicts \eqref{def_P_cont} and therefore this situation can not occur.
		
		Similarly, the case where $Q$ acts on one of the factors $W_1$ or $W_2$ in a trivial way produces a projection $\mathcal P$ that satisfies 
		$$(\dag\mathcal P)^2=-\frac{3}{2}\,\dag\mathcal P+10\,\mathcal P.$$
		This again is a contradiction to \eqref{def_P_cont}. Thus $Q\subset sp(2)$ must act irreducibly on $W,$ and therefore its action is given by  \eqref{rep_W} with respect to some appropriate $Sp(2)$-adapted frame $\big($because of the uniqueness of the irreducible representation in dimension $4\big).$

		Since, by \eqref{II_[pp]}, the projection  $\mathcal P$ is $Sp(1)$-equivariant, so is the tensor $S=\mathcal K^{-1}(K)$. It is well known that, with respect to the irreducible action of $Sp(1)$ in dimension $4$, there exists exactly one (up to scaling) equivariant totally symmetric 4-tensor. Therefore, $S$  must coincide with $\hat S$ from \eqref{new-dis}. This completes the proof.
	\end{proof}
	
	\begin{rmrk}\label{sp2_T_K_eigenspaces} The Lie algebra $sp(2)$ can be decomposed as $sp(2)=sp(1)_{ir}\oplus sp(1)_{ir}^\perp$.
		By Theorem~\ref{sp2-characterization}--(I), $sp(1)_{ir}$ and $sp(1)_{ir}^\perp$ are the eigenspaces  of $\mathcal T_K$ corresponding to the eigenvalues  $\frac{7}{2}$ and $-\frac{3}{2}$, respectively.
	\end{rmrk}

	\begin{rmrk}\label{sp2sp1-orbit}
		One immediate consequence of Theorem~\ref{sp2-characterization} is that the $Sp(2)$-orbit $\mathcal C$  is preserved under the action of the larger group $Sp(2)Sp(1)$ (cf. Remark~\ref{rmrk_spnsp1}), i.e., it is also a $Sp(2)Sp(1)$-orbit.
	\end{rmrk}
	
	\begin{cor}\label{cor_S_prop} \label{characterization_in_coordinates} Let $\{e_\alpha\}$ be any $Sp(2)$-adapted basis of $W$ and let $S\in\mathfrak S^{Sp(2)}$,  $S_{\alpha\beta\gamma\delta}=S(e_\alpha,e_\beta,e_\gamma,e_\delta),$ be any element. We have that $\mathcal K(S)\in\mathcal C$ if and only if the following two conditions are satisfied:

		\vspace{0.1cm}

		(I)  $\ $ $\pi^{\sigma\tau}\pi^{\mu\nu}S_{\sigma\mu\alpha\beta}S_{\tau\nu\gamma\delta}\ =\ \ 2\,S_{\alpha\beta\gamma\delta}\ +\ \frac{21}{8}\big(\pi_{\alpha\gamma}\pi_{\beta\delta}+\pi_{\alpha\delta}\pi_{\beta\gamma}\big)$
		
		\vspace{0.1cm}
		
		(II)  $\ $ $\sum_{\alpha,\beta,\gamma,\delta}\big(\pi^{\sigma\tau}S_{\sigma\alpha\beta\gamma}S_{\tau\delta\mu\nu}+ \frac{3}{4}S_{\alpha\beta\gamma\mu}\pi_{\nu\delta}+ \frac{3}{4}S_{\alpha\beta\gamma\nu}\pi_{\mu\delta}\big)\,x^\alpha x^\beta x^\gamma x^\delta\ =\ 0\\$
		for all $x^\ind{1},x^\ind{2},x^\ind{3},x^\ind{4}\in \mathbb C.$
		
		\vspace{0.2cm}
		
		Furthermore, if $K=\mathcal K(S)$ and $\mathcal P$ is defined by \eqref{def_P}, then equation (I) is equivalent to \eqref{I_p2} and (II) is equivalent to \eqref{II_[pp]}.
		
	\end{cor}
	
	\begin{rmrk}
		Formula (II) of Corollary~\ref{cor_S_prop} can also be expressed by stating that the symmetrization of the expression in the parentheses over the indices $\alpha$, $\beta$, $\gamma$ and $\delta$ is zero.
	\end{rmrk}

	\begin{proof}[Proof of the corollary] Identifying the complexification of $sp(2)$ with $\text{S}^2\,{W^*}$, we consider the basis  $$\{\$_{\alpha\beta}\ |\ 1\le\alpha\le\beta\le 4\}$$ form the proof of Lemma~\ref{T_K}. 
		
		If $K=\mathcal K(S)$, then 
		$
		\mathcal T_K(\$_{\alpha\beta})\ =\ S_{\alpha\beta\sigma\tau}\pi^{\sigma\gamma}\pi^{\tau\delta}\,\$_{\gamma\delta}
		$ and 
		the Lie bracket in $sp(2)$ is given by the second line in \eqref{Lie_br_sp2}. We calculate
		\begin{multline*}
			\Big[\mathcal T_K(\$_{\mu\nu}),\mathcal T_K(\$_{\alpha\beta})\Big]\ -\ \mathcal T_K\Big[\mathcal T_K( \$_{\mu\nu}),\mathcal \$_{\alpha\beta}\Big]\ -\ \frac{3}{2}\,\Big(\mathcal T_K\big[\$_{\mu\nu},\$_{\alpha\beta}\big]-\big[\$_{\mu\nu},\mathcal T_K(\$_{\alpha\beta})\big]\Big)\\
			=\ \Big\{\pi^{\sigma\tau}\big(S_{\sigma\mu\nu\alpha}S_{\tau\beta\gamma\delta}+S_{\sigma\mu\nu\beta}S_{\tau\alpha\gamma\delta}+S_{\sigma\mu\nu\gamma}S_{\tau\alpha\beta\delta}+S_{\sigma\mu\nu\delta}S_{\tau\alpha\beta\gamma}\big) 
			- \frac{3}{4}\big(S_{\alpha\beta\gamma\mu}\pi_{\nu\delta}+ S_{\alpha\beta\gamma\nu}\pi_{\mu\delta}\\
			+S_{\beta\gamma\delta\mu}\pi_{\nu\alpha}+ S_{\beta\gamma\delta\nu}\pi_{\mu\alpha}
			+S_{\gamma\delta\alpha\mu}\pi_{\nu\beta}+ S_{\gamma\delta\alpha\nu}\pi_{\mu\beta}
			+S_{\delta\alpha\beta\mu}\pi_{\nu\gamma}+ S_{\delta\alpha\beta\nu}\pi_{\mu\gamma}\big)\Big\}\,\pi^{\gamma\xi}\pi^{\delta\eta}\,\$_{\xi\eta},
		\end{multline*}
		which implies the equivalence between condition (II) of this proposition and condition (II) of Theorem~\ref{sp2-characterization}. 
		By a similar calculation, (I) here is equivalent to (I) in Theorem~\ref{sp2-characterization}.
		
	\end{proof}
	
	\begin{lemma} An element $S\in \mathfrak S^{Sp(2)}$ satisfies equation (II) of Proposition~\ref{characterization_in_coordinates} iff 
		\begin{equation}\label{sredni_indexi}
			\sum_{\alpha,\beta,\gamma,\delta}\big(\pi^{\sigma\tau}S_{\sigma\mu\alpha\beta}S_{\tau\nu\gamma\delta}- \frac{1}{4}S_{\alpha\beta\gamma\delta}\pi_{\mu\nu}+ \frac{1}{2}\pi_{\alpha\mu}S_{\nu\beta\gamma\delta}- \frac{1}{2}\pi_{\alpha\nu}S_{\mu\beta\gamma\delta}\big)\,x^\alpha x^\beta x^\gamma x^\delta\ =\ 0
		\end{equation}
		for all $x^\ind{1},x^\ind{2},x^\ind{3},x^\ind{4}\in \mathbb C$.
	\end{lemma}
	
	\begin{proof} By symmetrizing the expressions $\pi^{\sigma\tau}S_{\sigma\mu\alpha\beta}S_{\tau\nu\gamma\delta}$  and $\pi^{\sigma\tau}S_{\sigma\alpha\beta\gamma}S_{\tau\delta\mu\nu}$ over the indices $\alpha,\beta,\gamma,\delta$, we obtain, respectively, two 6-tensors:  $B_{\mu\nu|\alpha\beta\gamma\delta}\in \Lambda^2\,W^\ast\otimes\text{S}^4\,W^\ast$; and  $D_{\alpha\beta\gamma\delta|\mu\nu}\in \text{S}^4\,W^\ast\otimes\text{S}^2\,W^\ast$,
		\begin{equation*}
			\begin{aligned}
				B_{\mu\nu|\alpha\beta\gamma\delta}&=\pi^{\sigma\tau}\big(S_{\sigma\mu\alpha\beta}S_{\tau\nu\gamma\delta}+S_{\sigma\mu\gamma\delta}S_{\tau\nu\alpha\beta}+S_{\sigma\mu\alpha\gamma}S_{\tau\nu\beta\delta}+S_{\sigma\mu\beta\delta}S_{\tau\nu\alpha\gamma}\\
				&\hspace{8cm}+S_{\sigma\mu\alpha\delta}S_{\tau\nu\beta\gamma}+S_{\sigma\mu\beta\gamma}S_{\tau\nu\alpha\delta}\big),\\
				D_{\alpha\beta\gamma\delta|\mu\nu}&=\pi^{\sigma\tau}\big(S_{\sigma\alpha\beta\gamma}S_{\tau\delta\mu\nu}+S_{\sigma\alpha\beta\delta}S_{\tau\gamma\mu\nu}+S_{\sigma\alpha\gamma\delta}S_{\tau\beta\mu\nu}+S_{\sigma\beta\gamma\delta}S_{\tau\alpha\mu\nu}\big).
			\end{aligned}
		\end{equation*}
		It is a simple calculation to verify that
		\begin{equation*}
			\begin{aligned}
				&B_{\mu\nu|\alpha\beta\gamma\delta}=D_{\mu\alpha\beta\gamma|\delta\nu}-D_{\nu\alpha\beta\gamma|\delta\mu}\\
				&3D_{\alpha\beta\gamma\delta|\mu\nu}=B_{\alpha\mu|\beta\gamma\delta\nu}+B_{\beta\mu|\alpha\gamma\delta\nu}+B_{\gamma\mu|\alpha\beta\delta\nu}+B_{\delta\mu|\alpha\beta\gamma\nu}+2B_{\mu\nu|\alpha\beta\gamma\delta}
			\end{aligned}
		\end{equation*}
		
		\vspace{0.2cm}
		If we assume that the equation (II) of Proposition~\ref{characterization_in_coordinates} is satisfied, then
		\begin{equation*}
			\begin{aligned}
				D_{\alpha\beta\gamma\delta|\mu\nu}
				&= -\frac{3}{4}\Big(S_{\alpha\beta\gamma\mu}\pi_{\nu\delta}
				+S_{\alpha\beta\delta\mu}\pi_{\nu\gamma}
				+S_{\alpha\gamma\delta\mu}\pi_{\nu\beta}+S_{\beta\gamma\delta\mu}\pi_{\nu\alpha}
				\\
				&\hspace{2cm}+S_{\alpha\beta\gamma\nu}\pi_{\mu\delta}+S_{\alpha\beta\delta\mu}\pi_{\nu\gamma}+S_{\alpha\gamma\delta\mu}\pi_{\nu\beta}+S_{\beta\gamma\delta\mu}\pi_{\nu\alpha}\Big),
			\end{aligned}
		\end{equation*}
		and therefore
		\begin{equation*}
			\begin{aligned}
				6\pi^{\sigma\tau}S_{\sigma\mu\alpha\beta}S_{\tau\nu\gamma\delta}\,x^\alpha x^\beta x^\gamma x^\delta\ &= \ B_{\mu\nu|\alpha\beta\gamma\delta}\,x^\alpha x^\beta x^\gamma x^\delta\ =\ \big(D_{\mu\alpha\beta\gamma|\delta\nu}-D_{\nu\alpha\beta\gamma|\delta\mu}\big)\,x^\alpha x^\beta x^\gamma x^\delta\ \\[0.2cm]
				&=\big(\frac{3}{2}\,S_{\alpha\beta\gamma\delta}\pi_{\mu\nu}-3\pi_{\alpha\mu}S_{\nu\beta\gamma\delta}+3\pi_{\alpha\nu}S_{\mu\beta\gamma\delta}\big)\, x^\alpha x^\beta x^\gamma x^\delta
			\end{aligned}
		\end{equation*}
		
		For the converse, suppose that \eqref{sredni_indexi} is satisfied. Then, 
		\begin{equation*}
			\begin{aligned}
				B_{\mu\nu|\alpha\beta\gamma\delta}
				&= \frac{3}{2}\pi_{\mu\nu}S_{\alpha\beta\gamma\delta}+\frac{3}{4}\Big(\pi_{\mu\alpha}S_{\nu\beta\gamma\delta}
				+\pi_{\mu\beta}S_{\nu\alpha\gamma\delta}
				+\pi_{\mu\gamma}S_{\nu\alpha\beta\delta}+\pi_{\mu\delta}S_{\nu\alpha\beta\gamma}
				\\
				&\hspace{5cm}-\pi_{\nu\alpha}S_{\mu\beta\gamma\delta}
				-\pi_{\nu\beta}S_{\mu\alpha\gamma\delta}
				-\pi_{\nu\gamma}S_{\mu\alpha\beta\delta}-\pi_{\nu\delta}S_{\mu\alpha\beta\gamma}\Big)
			\end{aligned}
		\end{equation*}
		and we obtain
		\begin{equation*}
			\begin{aligned}
				&12\pi^{\sigma\tau}S_{\sigma\mu\alpha\beta}S_{\tau\nu\gamma\delta}\,x^\alpha x^\beta x^\gamma x^\delta\ = \ 3D_{\alpha\beta\gamma\delta|\mu\nu|}\,x^\alpha x^\beta x^\gamma x^\delta\\[0.2cm] 
				&\qquad=\ \big(4B_{\alpha\mu|\beta\gamma|\delta\nu}+2B_{\mu\nu|\alpha\beta\gamma\delta}\big)\,x^\alpha x^\beta x^\gamma x^\delta\
				=9\big(\pi_{\alpha\mu}S_{\nu\beta\gamma\delta}+\pi_{\alpha\nu}S_{\mu\beta\gamma\delta}\big)\, x^\alpha x^\beta x^\gamma x^\delta.
			\end{aligned}
		\end{equation*}
	\end{proof}

	\subsection{Stabilizer group}\label{isotropy}  Let  $G\subset Sp(2)$ be the stabilizer group of the element $\hat S\in\mathfrak S^{Sp(2)}$, given by \eqref{new-dis} (the cubic discriminant).  We claim that $G=Sp(1)_{\text{ir}}$. Indeed, by construction, $Sp(1)_{\text{ir}}$ is contained in $G$. On the other hand if $A$ is any element of $G$, then the adjoint action $Ad(A)$ of $A$ on $sp(2)$ commutes with ${\mathcal T}_{\hat K}$, $\hat K=\mathcal K(\hat S)$ \big(cf.~Lemma~\ref{1-1 HKC} and ~\eqref{def_K_1}\big). Therefore, $Ad(A)$  preserves the $\frac{7}{2}-$eigenspace of ${\mathcal T}_{\hat K}$, $sp(1)_{\text{ir}}\subset sp(2)$.  In other words, the restriction of $Ad(A)$ to $sp(1)_{\text{ir}}$ is an automorphism of $sp(1)_{\text{ir}}$. It is well known, that all automorphisms of $sp(1)$ are inner automorphisms, i.e., they are given by conjugation with an appropriate element of the Lie group $Sp(1)$ (cf. ~\cite{FH}, ~\cite{H}). 
	Therefore, there exists  $g\in Sp(1)_{\text{ir}}$  so that  
	\begin{equation*}
		A\cdot X\cdot A^{-1}=g\cdot X\cdot g^{-1}\qquad \forall X\in sp(1)_{\text{ir}}.
	\end{equation*}
	It follows that the conjugation with $B\overset{def}{=}g^{-1}\cdot A$ is trivial on $sp(1)_{\text{ir}}$, and therefore, by the Schur's lemma, $B$ is a scalar multiple of the identity, i.e., $A=\pm g$. Hence $A\in Sp(1)_{\text{ir}}$. In fact, we have shown the following: 
	
	\begin{prop} If $V$ is any 8-dimensional vector space equipped with a hyper-Hermitian structure, then the respective orbit $\mathcal C\subset \mathfrak R^{Sp(2)}$ is a 7-dimensional manifold, diffeomorphic to the homogeneous space
		$
		{Sp(2)}/{Sp(1)_{\text{ir}}}.
		$
		
	\end{prop}
	
	Alternatively, if we denote by $SO(4)_{ir}$ the product group $Sp(1)_{ir}Sp(1)\subset GL(V)$, then $SO(4)_{ir}$ is the stabilizer in $Sp(2)Sp(1)$ of $\hat K\in\mathcal C$,  and we have the representation 
	\begin{equation}\label{sp2sp1-so4}
		\mathcal C=\frac{Sp(2)Sp(1)}{SO(4)_{\text{ir}}}.
	\end{equation}
	
	\subsection{Irreducible representations of SO(4)}\label{irr_rep_so4} Let $E$ and $H$ be the standard irreducible representations on $\mathbb C^2$ of the first and second $Sp(1)$ factors, respectively,  in the decomposition
	\begin{equation}\label{so4-dec}
		SO(4)=\frac{Sp(1)\times Sp(1)}{\mathbb Z_2}.
	\end{equation} 
	It is well known (see for example \cite{Sal2} or \cite{Sw}) that the irreducible complex representations of $SO(4)$ are given, up to an isomorphism, by the formula
	\begin{equation}\label{rep-so4}
		\text{S}^k\,E\otimes \text{S}^l\,H,
	\end{equation}
	where $k$ and $l$ are non-negative integers satisfying $k\equiv l\ mod\ 2$. Each of these representations has an invariant real structure. The dimension of \eqref{rep-so4} is $(k+1)(l+1)$. 
	
	It follows that $SO(4)$ has exactly two irreducible representations on $\mathbb R^8$.  These are the real parts of the complex representations 
	$
	\text{S}^3 E\,\otimes\nobreak\,H$ and $E\,\otimes\, \text{S}^3H$. Clearly, in both cases, exactly one of the two $Sp(1)$ factors in the decomposition \eqref{so4-dec} acts irreducibly, while the other acts reducibly. Each of these two 8-dimensional representations  can be obtained from the other by a composition with a group authomorphism of $SO(4)$ that swaps the $Sp(1)$ factors in \eqref{so4-dec}.
	
	\subsection{The tangent space $T_K\mathcal C\subset \mathfrak R^{Sp(2)}$}
	Since the group $Sp(2)$ acts transitively on $\mathcal C\subset \mathfrak R^{Sp(n)}$ (cf. Definition~\ref{dfn_cubic_dis}), we can determine the tangent space at any point $K\in \mathcal C$ by examining the action of the Lie algebra $sp(2)$. The tangent space $T_K\mathcal C$ is composed of all elements $L\in\mathfrak R^{Sp(2)}$ that can be expressed as a Lie derivative of $K$. Specifically, for any $L\in T_K\mathcal C$,  there must exist some element $U\in sp(2),$ so that
	\begin{equation}\label{LinTK_U}
		L(x,y,z,w)\ =\ K(Ux,y,z,w)+K(x,Uy,z,w)+K(x,y,Uz,w)+K(x,y,z,Uw)
	\end{equation}
	for all vectors $x,y,z,w \in V.$

	\begin{lemma}\label{tangent_space_C}  Let $L$ be an element of the tangent space $T_K\mathcal C\subset \mathfrak R^{Sp(2)}$. We can define an endomorphism $H\in End(V)$ using the following formula:
		\begin{equation}\label{LinTK_H}
			\begin{aligned}
				g(Hx,y)=\frac{1}{120}\sum_{a,b,c=1}^8&\Big(L(x,h_a,h_b,h_c)\,K(y,h_a,h_b,h_c)\\
				& -\ L(y,h_a,h_b,h_c)\,K(x,h_a,h_b,h_c)\Big).
			\end{aligned}
		\end{equation}
		This holds for any $x,y\in V$ and any g-orthonormal basis $\{h_a\}_{a=1}^8$ of $V.$  Based on this definition, we have three key properties:  
		\par (i) $H\in sp(2)$.
		\par (ii) $\mathcal T_K(H)=-\frac{3}{2}H$. 
		\par (iii) $
		L(x,y,z,w)\ =\ K(Hx,y,z,w)+K(x,Hy,z,w)+K(x,y,Hz,w)+K(x,y,z,Hw)$
		for all $x,y,z,w \in V.$
		
		Furthermore, if  the expression for $L$ in \eqref{LinTK_U} holds for any other endomorphism  $U\in sp(2)$, then $H$ can be uniquely determined in terms of $U$ by the formula: $$H=\frac{1}{5}\Big(\frac{7}{2}U-\mathcal T_K(U)\Big).$$
	\end{lemma} 
	
	\begin{proof}  By a straightforward calculation, using Corollary~\ref{cor_S_prop}--(I), and Lemma~\ref{1-1 HKC}, one verifies that for each $K\in \mathcal C$, the following identities are satisfied: 
		\begin{equation}\label{contr_KxK_1}
			\begin{aligned}
				\sum_{a,b=1}^8 K(x,y,h_a,h_b)&K(z,w,h_a,h_b)\\
				&\ =\ 4K(x,y,z,w)\ +\ \frac{21}{8}\Big(g(x,z)g(y,w)-g(x,w)g(y,z)\Big)\\
				&\hspace{2cm} +\ \frac{21}{8}\sum_{s=1}^3\Big(g(I_sx,z)g(I_sy,w)-g(I_sx,w)g(I_sy,z)\Big)
			\end{aligned}
		\end{equation}
		and
		\begin{equation}\label{contr_KxK_2}
			\begin{aligned}
				\sum_{a,b=1}^8 K(x,h_a,h_b,y)&K(z,h_a,h_b,w)\ =\ 2K(x,z,y,w)\ +\ \frac{21}{8}\,g(x,z)g(y,w)\\
				&\qquad+\frac{21}{16}\,g(x,w)g(y,z)\ -\ \frac{21}{16}\sum_{s=1}^3g(I_sx,w)g(I_sy,z).
			\end{aligned}
		\end{equation}
		
		Let us now assume that \eqref{LinTK_U} is satisfied for some  $U\in sp(2)$. Using \eqref{contr_KxK_1} and \eqref{contr_KxK_2}, we calculate
		\begin{equation*}
			\begin{aligned}
				&\sum_{a,b,c=1}^8 L(x,h_a,h_b,h_c)\,K(y,h_a,h_b,h_c)\ =\ \sum_{a,b,c=1}^8\Big(K(Ux,h_a,h_b,h_c)\\
				&\hspace{2cm}+K(x,Uh_a,h_b,h_c)+K(x,h_a,Uh_b,h_c)+K(x,h_a,h_b,Uh_c)\Big)\,K(y,h_a,h_b,h_c)\\
				&\hspace{2cm}=\ -\ 4\sum_{a=1}^8\Big(K(x,h_a,y,Uh_a)+K(x,y,h_a,Uh_a)\Big)\ +\ 42g(Ux,y).
			\end{aligned}
		\end{equation*}
		It follows (by invoking \eqref{LinTK_H} and \eqref{def_K_0}) that 
		\begin{equation*}
			\begin{aligned}
				g(Hx,y)\ &=\ -\ \frac{1}{30}\sum_{a=1}^8\Big(K(x,h_a,y,Uh_a)-K(y,h_a,x,Uh_a)\\
				&\hspace{4cm}+\ 2\,K(x,y,h_a,Uh_a)\Big)\ +\ \frac{7}{10}g(Ux,y)\\
				&=\ -\frac{1}{10}\sum_{a=1}^8\, K(x,y,h_a,Uh_a)\ +\ \frac{7}{10}g(Ux,y)\ =\ \frac{1}{5}\,g\,\Big(\frac{7}{2}Ux-\mathcal T_K(U)x,\,y\Big).
			\end{aligned}
		\end{equation*}
		The rest is obvious.
	\end{proof}

	\section{$SO(4)_{ir}$ as a structure group}
	
	For any given irreducible representations  of $SO(4)$ on $\mathbb R^8$, we know there exists an $SO(4)$-equivariant inner product and an orientation. This representation corresponds to a subgroup of $SO(8)$ which we denote by $SO(4)_{ir }$.  A unique copy of the group $Sp(2)Sp(1)$ always exists within $SO(8)$ that contains $SO(4)_{ir }.$ To see how this works, we can follow the construction:
	\begin{itemize}
		\item[1.] Choose one of the two $Sp(1)$ factors from the decomposition \eqref{so4-dec} that acts reducibly on $\mathbb R^8$ (cf. Section~\ref{irr_rep_so4}).
		\item[2.] Find the centralizer of this chosen $Sp(1)$ factor within $SO(8).$ This centralizer turns out to be a subgroup isomorphic to $Sp(2).$
		\item[3.] The desired copy of $Sp(2)Sp(1)$ is the product of this $Sp(2)$ subgroup and the $Sp(1)$ factor from step 1.
	\end{itemize}
	
	In this section we shall consider oriented Riemannian 8-manifolds  admitting a reduction of the frame bundle to $SO(4)_{ir }\subset SO(8)$. By the above argument, such manifolds necessarily admit a (uniquely determined) reduction of the frame bundle to $Sp(2)Sp(1)\subset SO(8)$, i.e., they are always endowed with a canonical almost quaternion-Hermitian structure.
	
	\subsection{Almost quaternion-Hermitian manifolds}\label{almostQH} An almost quaternion-Hermitian 8-manifold is an 8-dimensional manifold equipped with a Riemannian metric $g$ and a rank-3 subbundle $\mathcal Q\subset End(TM)$ locally spanned by three endomorphisms $I_1,I_2,I_3$ that satisfy the identities 
	\begin{equation}\label{quat_id}
		I_1\cdot I_2=-I_2\cdot I_1=I_3;\qquad I_s\cdot I_s=-\text{Id};\qquad\text{and}\qquad g(I_sx,I_sy)=g(x,y)
	\end{equation}
	for all $x,y\in TM$ and $s=1,2,3.$  For a more in-depth discussion on this topic, please refer to \cite{Bes} or \cite{Sal2}.

	Consider the local 2-forms $\omega_s,\ s=1,2,3$, defined by $\omega_s(x,y)=g(I_sx,y)$, $\forall x,y\in TM$. The fundamental  4-form 
	\begin{equation}\label{def_Omega}
		\Omega=\sum_s\omega_s\wedge\omega_s
	\end{equation}
	is known to be independent of the choice of a local triple \eqref{quat_id}, and it is defined globally on $M$. If  $\Omega$ is parallel with respect to the Levi-Civita connection of $g$, then  $(M,g,\mathcal Q)$ is said to be a quaternion-K\"ahler manifold.
	If $\Omega$ is closed but not parallel (this can happen only in dimension 8, cf. \cite{Sw,CMS}), the structure  is called  harmonic. 
	
	Any almost quaternion-Hermitian 8-manifold $M$ admits a reduction of its linear frame bundle to a subgroup $Sp(2)Sp(1)\subset SO(8)$. This means there is always an associated principal $Sp(2)Sp(1)$-bundle $P(M)\rightarrow M$ which is a subbundle of the linear frame bundle of $M$.
	Let us fix (once and for all)  a triple of (constant) endomorphisms $J_1,J_2,J_3$ on $V=\mathbb R^8$ that are compatible with the standard inner product there and satisfy \eqref{quat_iden}. With respect to $J_1$, the complexification $V^{\mathbb C}$ of $V$ splits into a direct sum of  $i$ and $-i$ eigenspaces, \begin{equation}\label{v=w+wb}
		V^{\mathbb C}=W+\overline{W}.
	\end{equation} 
	
	We shall use also a fixed $Sp(2)$-adapted (cf. \eqref{constants}) basis $\{e_{\alpha}\}$ of $W$ (and respectively of $V^{\mathbb C}$ with the conventions introduced in Section~\ref{index_conventions}). The tautological 1-form $\theta$ on $P(M)$ takes values in $V^{\mathbb C}$. The composition of $\theta$ with $e^{\alpha}\in W^\ast$ and $e^{\bar\alpha}\in\overline{W}^\ast$ produces a set of globally defined 1-forms $\theta^{\alpha},\theta^{\bar\alpha}$ on $P(M)$.

	If $M$ is a quaternion-K\"ahler manifold, then the Levi-Civita connection $\nabla$ of $g$ is represented by a 1-form on $P(M)$ with values in the Lie algebra $sp(2)+sp(1)$. Using the identification (cf. \eqref{Spn_X})
	\begin{equation}
		sp(2)=\Big\{(X_{\alpha\beta})_{4\times 4}\ \Big|\ X_{\alpha\beta}=X_{\beta\alpha},\ (\mathfrak jX)_{\alpha\beta}=X_{\alpha\beta}\Big\}
	\end{equation}
	and 
	\begin{equation}\label{sp1_rep}
		sp(1)=\Big\{\sum_{s=1}^3a^s\,J_s\ \Big|\ a^s\in\mathbb R\Big\},
	\end{equation}
	we can represent  $\nabla$ on $P(M)$ by
	\begin{equation}\label{Levi-Civita}
		\Big((\Gamma_{\alpha\beta})_{4\times 4},\,\frac{1}{2}\sum_s \phi^s J_s\Big)
	\end{equation}
	(the first entry represents the $sp(2)$-component of the connection and the second---the $sp(1)$-component), where $\phi^s$ are real-valued and $\Gamma_{\alpha\beta}$ are complex-valued 1-forms that satisfy $\Gamma_{\alpha\beta}=\Gamma_{\beta\alpha}$ and $ (\mathfrak j\Gamma)_{\alpha\beta}\footnotemark=\Gamma_{\alpha\beta}$\footnotetext{For tensor valued 1-forms we follow, as well, the conventions from Section~\ref{index_conventions}, e.g., we have that: $(\mathfrak j\Gamma)_{\alpha\beta}=\pi^{\bar\sigma}_{\dt\alpha}\pi^{\bar\tau}_{\dt\beta}\Gamma_{\bar\sigma\bar\tau}$; and $\Gamma_{\bar\sigma\bar\tau}=\overline{\Gamma_{\sigma\tau}}.$}. According to our conventions and the chosen (in Section~\ref{preliminaries}) matrix representation for the Lie algebra $sp(2) + sp(1)$ on $V$, any vector $\theta^\alpha\,e_\alpha + \theta^{\bar\alpha}\,e_{\bar\alpha}\in V$ is mapped by \eqref{Levi-Civita} to a vector $\eta^\alpha\,e_\alpha + \eta^{\bar\alpha}\,e_{\bar\alpha}\in V$, where
	\begin{equation*}
		\eta^\alpha\ =\ \pi^{\alpha\sigma}\,\Gamma_{\sigma\beta}\,\theta^{\beta}\ +\ \frac{i}{2}\,\phi^1\,\theta^\alpha\ -\ \frac{1}{2}\,\pi^\alpha_{\,.\, \bar\beta}\,\big(\phi^2+i\phi^3\big)\,\theta^{\bar\beta}
	\end{equation*}
	and, of course, $\eta^{\bar\alpha}=\overline{\eta^\alpha}$. Since, by definition, $\nabla$ is torsion free, we have the structure equations 
	\begin{equation}\label{str_eq_LC}
		d\theta^{\alpha}\ =\ -\pi^{\alpha\sigma}\,\Gamma_{\sigma\beta}\wedge\theta^{\beta}\ -\ \frac{i}{2}\,\phi^1\wedge\theta^\alpha\ +\ \frac{1}{2}\,\pi^\alpha_{\,.\, \bar\beta}\,\big(\phi^2+i\phi^3\big)\wedge\theta^{\bar\beta}.
	\end{equation}

	The curvature tensor 
	\begin{equation}
		R(,)=\Big[\nabla,\nabla\Big]-\nabla_{[,]}
	\end{equation}
	of any quaternion-K\"ahler manifold $M$ splits canonically into a sum (cf. \cite{Ale2} or \cite{Sal}, Theorem. 3.1)
	\begin{equation}\label{gen_cur_split}
		R=R'+\frac{Scal}{64} R_0
	\end{equation}
	where $R'$ is a hyper-K\"ahler curvature type tensor (pointwise on $M$), $Scal$ is the scalar curvature of $g$ (which is a constant) and $R_0$ is a tensor field given by
	\begin{equation}\label{proj-cur}
		\begin{aligned}
			4g\Big(R_0(x,y)z,&w\Big)\ =\  g(x,w)g(y,z)-g(x,z)g(y,w)\\
			&+\ \sum_{s=1}^3\Big(-2\omega_s(x,y)\omega_s(z,w)+\omega_s(x,z)\omega_s(w,y)+\omega_s(x,w)\omega_s(y,z)\Big),
		\end{aligned}
	\end{equation} 
	for all $x,y,z,w\in TM$  (see for example ~\cite{Bes} 14.44).
	Notice that \eqref{proj-cur} is precisely the formula for the curvature of the quaternionic projective space $\mathbb H\mathbb P^n$.  The decomposition \eqref{gen_cur_split} implies, in particular, that each quaternion-K\"ahler manifold is Einstein.
	
	\subsection{Cubic discriminants on 8-manifolds}\label{cub_dis_sec} Suppose that $(M,g,\mathcal Q)$ is an 8-dimensional  almost quaternion-Hermitian manifold.  The adjoint action of the structure group $Sp(2)Sp(1)$ on $sp(2)$ crates  a vector bundle, 
	\begin{equation}\label{Bsp2}
		Bsp(2)=P(M)\times_{Ad}sp(2).
	\end{equation}  
	$Bsp(2)$ is a subbundle of the endomorphism bundle $End(TM)$. Its fiber, $Bsp(2)_p$, over any point $p\in M$, consist of all skew-symmetric endomorphisms $A\in End(T_pM)$ that commute with every element in $\mathcal Q_p$. Consequently, each $Bsp(2)_p$ forms  a Lie algebra isomorphic to $sp(2)$. 

	We may consider the space $\mathfrak R^{Sp(2)}_p$ of all hyper-K\"ahler curvature type tensors at a point $p$, cf. \eqref{HK_curv_prop}.  The corresponding vector bundle over $M$ is denoted by $\mathfrak R^{Sp(2)}(M)\rightarrow M.$  For each $K\in \mathfrak R^{Sp(2)}_p$, we have the linear endomorphism $\mathcal T_K\in End(Bsp(2)_p)$ defined as in \eqref{def_K_2}. 
	We let $\mathcal C_p$ to be the respective orbit of cubic discriminants at $p$ \big(cf. Definition~\ref{dfn_cubic_dis} and Theorem~\ref{sp2-characterization}\big), i.e.,  the set of all $K\in \mathfrak R^{Sp(2)}_p$ for which the following two conditions are satisfied:
	\begin{equation}\label{def_I_cd}
		\big(2\mathcal T_K-7\,\text{Id}\big)\big(2\mathcal T_K+3\,\text{Id}\big)\ =\ 0;
	\end{equation}
	\begin{equation}\label{def_II_cd}
		\begin{aligned} 
			\ \text{and}\\
			&\big[\mathcal T_KA,\mathcal T_KB\big]-\mathcal T_K\big[\mathcal T_KA, B\big]\ =\ \frac{3}{2}\,\Big(\mathcal T_K\big[A,B\big]-\big[A,\mathcal T_KB\big]\Big)\\
			\text{for all }A,B\in &Bsp(2)_p.
		\end{aligned}
	\end{equation}
	
	By choosing a frame  at $p$ from the principal fiber bundle $P(M)$, we obtain a representation of the structure group $Sp(n)Sp(1)$ on $T_pM$ that preserves $\mathcal Q_p$.  By \eqref{sp2sp1-so4}, we have that
	\begin{equation*}
		\mathcal C_p\cong\frac{Sp(2)Sp(1)}{SO(4)_{ir}}.
	\end{equation*}
	The 7-dimensional orbits  $\mathcal C_p$, $p\in M$,  are the fibers of a certain locally trivial fiber bundle $\mathcal C(M)$ over $M$. The smooth sections of $\mathcal C(M)$ will be called cubic discriminants on~$M$.

	If $\sigma$ is any local section of the principal fiber bundle $P(M)$, the pullback of the tautological 1-form on $P(M)$ via $\sigma$ defines a local coframing 
	\begin{equation}\label{sigma_cof}
		\sigma^\ast\theta^\alpha,\  \sigma^\ast\theta^{\bar\alpha}
	\end{equation} on $M$. For any section $K$ of $\mathfrak R^{Sp(2)}(M)$, we define 
	\begin{equation}
		S_{\alpha\beta\gamma\delta} = K_{\alpha{\bar\sigma}\gamma{\bar\tau}}\, \pi^{\bar\sigma}_{\ \beta}\pi^{\bar\tau}_{\ \delta},
	\end{equation}
	where $K_{\alpha{\bar\sigma}\gamma{\bar\tau}}$ are the components of $K$ with respect to the coframing \eqref{sigma_cof}. By Lemma~\ref{1-1 HKC}, $S_{\alpha\beta\gamma\delta}$ is a totally symmetric tensor field and $S_{\alpha\beta\gamma\delta}=(\mathfrak j S)_{\alpha\beta\gamma\delta}$. The 4-tensor $K$ is a cubic discriminant iff around each point of $M$ there exist a section $\sigma$ of $P(M)$, defined on a open subset $U\subset M$, so that 
	\begin{equation}\label{cub-filed-symm}
		\begin{aligned}
			S_{\alpha\beta\gamma\delta}\,x^\alpha x^\beta x^\gamma x^\delta&\ =\ -18\, x^\ind{1} x^\ind{2} x^\ind{3} x^\ind{4}+4\sqrt{3}\,x^\ind{1}(x^\ind{4})^\pow{3}\\
			&
			-4\sqrt{3}\,(x^\ind{2})^\pow{3}x^\ind{3}+3(x^\ind{2})^\pow{2}(x^\ind{4})^\pow{2}-9(x^\ind{1})^\pow{2}(x^\ind{3})^\pow{2}
		\end{aligned}
	\end{equation}
	at each point of $U$ and for all $x^1,x^2,x^3,x^4\in \mathbb C$ (cf. Remark~\ref{rmk-cd}).

	Having a globally defined cubic discriminant on $M$ is the same as reducing the structure group of $M$ from $Sp(2)Sp(1)$ to $SO(4)_{ir }$. This means  that for every global section $K$ of $\mathcal C(M),$ we can create a unique principal $SO(4)_{ir}$-bundle, denoted as $\mathcal D(M)$, which is a subbundle of $P(M).$ This new bundle,
	\begin{equation}\label{DM}
		\mathcal D(M)\rightarrow M,
	\end{equation}
	is constructed so that its sections are precisely the sections $\sigma$ of $P(M)$ for which \eqref{cub-filed-symm} is satisfied.  The structure group $SO(4)_{ir}$ of $\mathcal D(M)$ can be understood as the stabilizer subgroup of $K$  within $Sp(2)Sp(1)$ at any given point of the manifold. The adjoint action of $SO(4)_{ir}$ on $so(4)_{ir}=sp(1)_{ir}\oplus sp(1)$ preserves the ideal $sp(1)_{ir}$ and determines  a vector bundle  
	\begin{equation}\label{bndl-sp1ir}
		Bsp(1)_{ir}\ =\ D(M)\times_{Ad}sp(1)_{ir} \ \subset Bsp(2)\subset End(TM).
	\end{equation}

	Drawing on the same method as in the discussion surrounding  \eqref{v=w+wb}, we can use the tautological 1-form $\theta$ on $D(M)$ to create a set of globally defined (complex valued) 1-forms $\theta^{\alpha},\theta^{\bar\alpha}$ on $\mathcal D(M)$. If $\sigma$ is any local sections of $\mathcal D(M)$ and 
	\begin{equation}\label{sigma_cof_D}
		\sigma^\ast\theta^\alpha,\  \sigma^\ast\theta^{\bar\alpha}
	\end{equation} is the corresponding local coframing on $M$, then $Bsp(1)_{ir}$  is locally the linear span of a triple $\mathcal E_1,$ $\mathcal E_2$, $\mathcal E_3$ of pointwise endomorphisms of $TM$ that are defined w.r.t. \eqref{sigma_cof_D} by 
	\begin{equation}\label{gen_E}
		\mathcal E_s\mapsto\pi^{\alpha\sigma}(\Upsilon_s)_{\sigma\beta}, \qquad s=1,2,3.
	\end{equation}
	where $\Upsilon_s$ are the matrices \eqref{upsilon}.  Unlike the local generators $I_s$ of $\mathcal Q$ \big(cf. \eqref{quat_id}\big), the endomorphisms $\mathcal E_s$ do not square to -Id (they are not almost complex structures on $M$). They possess, however, a number of other useful identities, e.g., the following three that come directly from Lemma~\ref{upsilon_lemma}---(2), (3) and (7).
	\begin{lemma}\label{E_s_id}
		We have the identities:
		\begin{equation}\label{1_E_s_id}
			g\big(\mathcal E_s,\mathcal E_t\big)\ =\ 
			\begin{cases}
				5,\ \text{if}\ s=t,\\
				0,\ \text{otherwise}.\footnotemark
			\end{cases}
		\end{equation} \footnotetext{The inner product $g$ is extended over the endomorphisms of $TM$ by the formula $g(\mathcal A,\mathcal B)=\frac{1}{2}\text{trace}(A^TB)$, where $\mathcal A$ and $\mathcal B$ have matrices $ A$ and $ B$ w.r.t. any $g$-orthonormal framing on $M$.}
		\begin{equation}\label{2_E_s_id}
			\big[\mathcal E_i,\mathcal E_j\big]\ =\ \mathcal E_k;\ \footnotemark
		\end{equation}\footnotetext{$(ijk)$ is any positive permutation of $(123)$.}
		\begin{equation}\label{3_E_s_id}
			\mathcal E_1\cdot \mathcal E_1+\mathcal E_2\cdot \mathcal E_2+\mathcal E_3\cdot \mathcal E_3\ =\ -\frac{15}{4}\,\text {Id}.
		\end{equation}
		
	\end{lemma}
	
	Furthermore, using Lemma~\ref{upsilon_lemma}---(5), we easily calculate that 
	\begin{equation}\label{formula_K_e_omega}
		\begin{aligned}
			K(x,y,z,w)\ =\ &\sum_s{g(\mathcal E_sx,y)g(\mathcal E_sz,w)}\ +\ \frac{3}{8}\Big(g(x,w)g(y,z)-g(x,z)g(y,w)\Big)\\
			&+\ \frac{3}{8}\sum_{s=1}^3\Big(\omega_s(x,z)\omega_s(w,y)+\omega_s(x,w)\omega_s(y,z)\Big)
		\end{aligned}
	\end{equation}
	for all $x,y,z,w\in TM$, where $\omega_s$ are as in \eqref{def_Omega} and $I_s$ are any local generators for $\mathcal Q$ that satisfy \eqref{quat_id}. In particular, if we set $\epsilon_s(x,y)=g(\mathcal E_sx,y)$, then  by an antisymmetrization of \eqref{formula_K_e_omega}, we obtain that
	\begin{equation}\label{e_s_to_Omega}
		\sum_s \epsilon_s\wedge\epsilon_s\ =\ -\ \frac{3}{4}\,\Omega,
	\end{equation}
	where $\Omega$ is the fundamental 4-form given by \eqref{def_Omega}.
	In fact, as we show below, just the two equations \eqref{e_s_to_Omega} and \eqref{2_E_s_id} alone are sufficient to determine a cubic discriminant on $M$.
	\begin{prop} Suppose that $(M,g,\mathcal Q)$ is an almost quaternion-Hermitian $8$-manifold.  Let $\mathcal E_s$, $s=1,2,3,$ be any globally defined  sections of  $Bsp(2)$ that satisfy \eqref{2_E_s_id}.  There exists a cubic discriminant $K$ on $M$, so that the induced bundle $Bsp(1)_{ir}$ $\big($defined by \eqref{bndl-sp1ir}$\big)$ coincides with the linear span of $\mathcal E_1$, $\mathcal E_2$ and $\mathcal E_3$, if and only if \eqref{e_s_to_Omega} is satisfied.
		
		In the case when \eqref{e_s_to_Omega} is satisfied, the respective cubic discriminant $K$ on $M$ is given by formula \eqref{formula_K_e_omega}. 
		
	\end{prop}
	
	\begin{proof}
		We have discussed already the necessity of \eqref{e_s_to_Omega}. To prove the converse, let us pick an arbitrary point $p\in M$ and let $\{h_a\}_{a=1}^8$ be any $g$-orthonormal basis for $T_pM$. Using \eqref{def_Omega}, we calculate
		\begin{equation*}
			\begin{aligned}
				\sum_{a=1}^{8}\Omega(x,y,h_a,I_1h_a)\ &=\ 2\Big(\omega_s(x,y)\omega_s(h_a,I_1h_a)\ +\ \omega_s(y,h_a)\omega_s(x,I_1h_a)\\
				&+\omega_s(h_a,x)\omega_s(y,I_1h_a)\Big)\ =\ 20\,\omega_1(x,y),\qquad\forall x,y\in T_pM.
			\end{aligned}
		\end{equation*}
		If instead of \eqref{def_Omega}, we use \eqref{e_s_to_Omega}, by performing a similar calculation, we obtain that
		\begin{equation*}
			\sum_{a=1}^{8}\Omega(x,y,h_a,I_1h_a)\ =\ -\frac{16}{3}\sum_{s=1}^3g(\mathcal E_sI_1\mathcal E_s\,x,y).
		\end{equation*}
		Since, by assumption, $\mathcal E_s$ are sections of $Bsp(2)$, the above calculation implies the validity of \eqref{3_E_s_id}.

		On the other hand, equation \eqref{2_E_s_id} shows that we are given a representation of the Lie algebra $sp(1)$ on the 8-space $T_pM$. The inner product on $End(T_pM)$, defined by 
		\begin{equation*}
			g(X,Y)=\frac{1}{2}\sum_{a,b=1}^8g(Xh_a,h_b)g(Yh_a,h_b), \qquad X,Y\in End(T_pM),
		\end{equation*} is apparently $sp(1)$-equivariant and therefore its restriction to the image of $sp(1)$ in $End(T_pM)$ must be proportional to the Killing form of $sp(1)$; that is, there exist a constant $l\in \mathbb R$ so that
		\begin{equation}\label{l_killing}
			g\big(\mathcal E_s,\mathcal E_t\big)\ =\ 
			\begin{cases}
				l,\ \text{if}\ s=t,\\
				0,\ \text{otherwise}.
			\end{cases}
		\end{equation} 
		By taking the trace of \eqref{3_E_s_id} and applying \eqref{l_killing}, we conclude that $l=5$ and thus \eqref{1_E_s_id} holds.
		
		Let us now consider an endomorphism $\mathcal P\in End(Bsp(2)_p)$, defined by
		\begin{equation*}
			\mathcal P(X)\ =\ \frac{1}{5}\sum_{s=1}^3\mathcal E_s\,g(\mathcal E_s,X),\qquad  X\in Bsp(2)_p.
		\end{equation*}
		It follows from \eqref{1_E_s_id} and \eqref{2_E_s_id} that
		\begin{equation}\label{pr_I_p12}
			\mathcal P^2=\mathcal P\qquad\text{and}\qquad \Big[\mathcal PX,\mathcal P Y\Big]=\mathcal P\Big[\mathcal PX,Y\Big],\quad \forall X,Y\in Bsp(2)_p.
		\end{equation} 
		
		The tensor $K$, defined by formula \eqref{formula_K_e_omega}, is obviously a hyper-K\"ahler curvature type tensor on $T_pM$. A straightforward calculation, using \eqref{formula_K_e_omega}, shows that
		\begin{equation}\label{pr_tk_p}
			\mathcal P{\ =\ }\frac{1}{5}\mathcal T_K+\frac{3}{10}\text {Id}.,
		\end{equation} 
		where $\mathcal T_K$ is defined as in \eqref{def_K_0}. By substituting \eqref{pr_tk_p} into \eqref{pr_I_p12} and expanding, we obtained 
		that $\mathcal T_K$ satisfies \eqref{def_I_cd} and \eqref{def_II_cd}, and thus $K$ is a cubic discriminant on $M$.
	\end{proof}

	\subsection{Two symmetric examples}\label{wolf}  We shall look at two examples of quaternion-K\"ahler manifolds---symmetric spaces of the form $G/SO(4)$--- that come with a built-in cubic discriminant. In these examples, $G$ is either the compact exceptional Lie group $G_2$ or its  noncompact real form, $G_{2(2)}$.
	
	Following \cite{FH}, the complexified Lie algebra $g_2\otimes\mathbb C$ ($g_2$ is the Lie algebra of the compact $G_2$) has a 2-dimensional Cartan subalgebra $h\subset g_2\otimes\mathbb C$ with two primitive roots
	\begin{equation*}
		\alpha_1=(1,0)\qquad \text{and}\qquad \alpha_2=(-\frac{3}{2},\frac{\sqrt{3}}{2}).
	\end{equation*} 
	The other positive roots of $g_2$ are given by
	\begin{equation*}
		\alpha_3=\alpha_1+\alpha_2,\ \alpha_4=2\alpha_1+\alpha_2,\ \alpha_5=3\alpha_1+\alpha_2,\ \alpha_6=3\alpha_1+2\alpha_2.
	\end{equation*} 
	There is a basis (cf. \cite{FH}, page 346) 
	\begin{equation*}
		H_1,H_2,X_1,X_2,X_3,X_4,X_5,X_6,Y_1,Y_2,Y_3,Y_4,Y_5,Y_6
	\end{equation*}
	of $g_2\otimes\mathbb C$, attached to the choice  of roots made above; the multiplication table for this basis is explicitly given in \cite{FH}, Table~22.1. 
	
	It is a well known fact that $g_2\otimes\mathbb C$ has exactly two real forms, up to isomorphism: the compact form $g_2$ and the split form $g_{2(2)}$. The compact form is associated with a specific conjugation map that sends $H_i$ to $-H_i$, $X_i$ to $-Y_i$ and $Y_i$ to $-X_i$. For our purposes, when dealing with the compact form,  the following basis will be useful:
	\begin{equation}\label{g_2_basis}
		\begin{aligned}
			\psi_1&=\frac{i}{2}H_1,\qquad&\psi_2&=\frac{1}{2}(X_1-Y_1),\qquad&\psi_3&=\frac{i}{2}(X_1+Y_1),\\
			\varphi_1&=\frac{i}{2}H_1+iH_2,\qquad&\varphi_2&=\frac{1}{2}(X_6-Y_6),\qquad&\varphi_3&=\frac{i}{2}(X_6+Y_6),
		\end{aligned}
	\end{equation}
	\begin{equation*}
		\begin{aligned}
			&e_1=\frac{\sqrt{3}}{\sqrt{2}}\,X_2,&\qquad &e_2=\frac{1}{\sqrt{2}}\,X_3,&\qquad &e_3=-\frac{\sqrt{3}}{\sqrt{2}}\,X_5,&\qquad &e_4=-\frac{1}{\sqrt{2}}\,X_4,&\\
			&e_{\bar 1}=-\frac{\sqrt{3}}{\sqrt{2}}\,Y_2,& &e_{\bar 2}=-\frac{1}{\sqrt{2}}\,Y_3,& &e_{\bar 3}=\frac{\sqrt{3}}{\sqrt{2}}\,Y_5,& &e_{\bar 4}=-\frac{1}{\sqrt{2}}\,Y_4.&
		\end{aligned}
	\end{equation*}

	With respect to the dual to \eqref{g_2_basis} basis $\psi^1,\psi^2,\psi^3,\varphi^1, \varphi^2, \varphi^3$, $e^\alpha,e^{\bar \alpha}$, the structure equations of  $G_2$ take the form:
	\begin{equation}\label{g2_dpsi}
		\begin{aligned}
			d\psi^1&\ =\ -\psi^2\wedge\psi^3\ +\ \pi^\sigma_{\,.\,\bar\beta}\,(\Upsilon_1)_{\alpha\sigma}\,e^{\alpha}\wedge e^{\bar\beta}\\
			d\psi^2&\ =\ -\psi^3\wedge\psi^1\ +\ \pi^\sigma_{\,.\,\bar\beta}\,(\Upsilon_2)_{\alpha\sigma}\,e^{\alpha}\wedge e^{\bar\beta}\\
			d\psi^3&\ =\ -\psi^1\wedge\psi^2\ +\ \pi^\sigma_{\,.\,\bar\beta}\,(\Upsilon_3)_{\alpha\sigma}\,e^{\alpha}\wedge e^{\bar\beta}\\
		\end{aligned}
	\end{equation}

	\begin{equation}\label{g2_dvarphi}
		\begin{aligned}
			d\varphi^1&\ =\ -\varphi^2\wedge\varphi^3\ -\ \frac{3i}{2}\Big(e^1\wedge e^{\bar 1}+e^2\wedge e^{\bar 2}+e^3\wedge e^{\bar 3}+e^4\wedge e^{\bar 4}\Big),\\
			d\varphi^2&\ =\ -\varphi^3\wedge\varphi^1\ -\ \frac{3}{2}\Big(e^1\wedge e^{3}+e^{2}\wedge e^{4}+e^{\bar 1}\wedge e^{\bar 3}+e^{\bar 2}\wedge e^{\bar 4}\Big),\\
			d\varphi^3&\ =\ -\varphi^1\wedge\varphi^2\ +\ \frac{3i}{2}\Big(e^1\wedge e^{3}+e^{2}\wedge e^{4}-e^{\bar 1}\wedge e^{\bar 3}-e^{\bar 2}\wedge e^{\bar 4}\Big),
		\end{aligned}
	\end{equation}

	\begin{equation}\label{dealpha}
		\begin{aligned}
			de^\alpha&\ =\ -\pi^{\alpha\sigma}\,(\Upsilon_s)_{\sigma\beta}\,\psi^s\wedge e^\beta\ -\ \frac{i}{2}\,\varphi^1\wedge e^\alpha\ +\ \frac{1}{2}\,\pi^\alpha_{.\, \bar\beta}\,\big(\varphi^2+i\varphi^3\big)\wedge e^{\bar\beta}, \\
			de^{\bar\alpha}&\ =\ -\pi^{\bar\alpha\bar\sigma}\,(\Upsilon_s)_{\bar\sigma\bar\beta}\,\psi^s\wedge e^{\bar\beta}\ +\ \frac{i}{2}\,\varphi^1\wedge e^{\bar\alpha}\ +\ \frac{1}{2}\,\pi^{\bar\alpha}_{.\, \beta}\,\big(\varphi^2-i\varphi^3\big)\wedge e^{\beta},
		\end{aligned}
	\end{equation}
	where $\Upsilon_s$ are the matrices from \eqref{upsilon_E}. A subalgebra of $g_2$ is formed by the real linear combinations of the basis vectors $\psi_1,\psi_2,\psi_3$, $\varphi_1, \varphi_2$ and $\varphi_3$. The corresponding closed subgroup of $G_2$ is isomorphic to $SO(4)$. We are interested in studying the homogeneous structure of the resulting quotient space, $G_2\big/ SO(4)$. 
	
	Consider the complex subspace $W\subset g_2\otimes \mathbb C$ with basis $e_1,e_2,e_3,e_4$. The conjugated  space $\overline{W}$ is generated by $e_{\bar 1},e_{\bar 2},e_{\bar 3},e_{\bar 4}$. If we denote by $V$ the real part of $W\oplus\overline {W}$ then $V$ can be identified with the tangent space at the origin in $G_2/SO(4)$. The $Ad(SO(4))$-invariant inner product $\sum \big(e^\alpha\otimes e^{\bar\alpha}+e^{\bar\alpha}\otimes e^{\alpha}\big)$ on $V$ defines a homogeneous  Riemannian metric tensor $g$ on $G_2/SO(4)$. The adjoint action of $\psi_1,$ $\psi_2$ and $\psi_3$ preserves $W$ and is given by the matrices $\pi^{\alpha\sigma}(\Upsilon_1)_{\sigma\beta}$, $\pi^{\alpha\sigma}(\Upsilon_2)_{\sigma\beta}$ and $\pi^{\alpha\sigma}(\Upsilon_3)_{\sigma\beta}$. The adjoint action of $\varphi_1$,  $\varphi_2$ and $\varphi_3$ commutes with the complex conjugation on $W\oplus\overline {W}$. If we denote by $J_1$, $J_2$ and $J_3$ the endomorphisms of $V$ that correspond to the adjoint action of $2\varphi_1, 2\varphi_2$ and $2\varphi_3$, then $W$ and $\overline{W}$ are the i and -i eigenspaces of $J_1$ and we have the identities \eqref{quat_iden}.
	The linear span of the endomorphisms $J_1$, $J_2$ and $J_3$ is a subspace in $End(V)$ that remains invariant under the adjoint action of $SO(4)$. Therefore, it determines a subbundle $\mathcal Q$ of the endomorphism bundle of $G_2/SO(4)$,  which is an almost quaternion-Hermitian structure. It is a well known fact that $(g,\mathcal Q)$ is in fact a quaternion-K\"ahler structure.
	
	We may regard $G_2\rightarrow G_2/SO(4)$ as a principal SO(4)-bundle associated with the canonical $Sp(2)Sp(1)$-bundle  (determined by $g$ and $\mathcal Q$). The 1-forms $e^\alpha$, $e^{\bar\alpha}$ represent the tautological ($V^{\mathbb C}$-valued) 1-form. The 1-forms $(\psi^s,\varphi^s)$ define a $so(4)$-valued connection on this principal bundle. The structure equations $\eqref{dealpha}$ imply that this connection is torsion free, i.e., it is the Levi-Civita connection $\nabla$ on $G_2/SO(4)$. Since $G_2/SO(4)$ is a symmetric space, the curvature $R$ of $\nabla$ is a $\nabla$-paralell tensor. For any $X,Y\in V$, $R(X,Y)$ is given by the projection onto $so(4)\subset g_2$ of $-[X,Y]$, i.e., we have
	\begin{equation}\label{R_sym}
		R=\sum_s\Big(d\psi^s|_{\Lambda^2V}\otimes ad(\psi_s)+d\varphi^s|_{\Lambda^2V}\otimes ad(\varphi_s)\Big).
	\end{equation}
	
	We can define a natural cubic discriminant on the manifold $G_2/SO(4)$. This discriminant is derived as a specific component of the manifold's curvature tensor, as will be shown in the next proposition.
	\begin{prop}\label{wolf_curvature} If $R=R'+\frac{Scal}{64}R_0$ is the splitting, in accordance with \eqref{gen_cur_split},  of the Riemannian curvature tensor of the symmetric space $G_2/SO(4)$, then $-R'$ is a parallel cubic discriminant.
	\end{prop}
	\begin{proof}We calculate using \eqref{R_sym},\eqref{g2_dpsi}, \eqref{g2_dvarphi} and Lemma~\eqref{upsilon_lemma}--\it{(5)}:
		\begin{equation*}
			\begin{aligned}
				g\Big(R(e_\alpha,e_{\bar\beta})e_\gamma,e_{\bar\delta}\Big)&=d\psi^s(e_\alpha,e_{\bar\beta})g\Big(\pi^{\mu\tau}\
				(\Upsilon_s)_{\tau\gamma}\,e_{\mu},e_{\bar\delta}\Big)+d\varphi_1(e_\alpha,e_{\bar\beta})g\Big(\frac{1}{2}J_1e_\gamma, e_{\bar\delta}\Big)\\
				&=-\pi^\sigma_{\,.\,\bar\beta}\,\pi^\tau_{\,.\,\bar\delta}\,\Big((\Upsilon_s)_{\alpha\sigma}\,(\Upsilon_s)_{\gamma\tau}\Big)+\ \frac{3}{4}\,g_{\alpha\bar\beta}\,g_{\gamma\bar\delta}\\
				&=-\pi^\sigma_{\,.\,\bar\beta}\pi^\tau_{\,.\,\bar\delta}\Big(\hat S_{\alpha\sigma\gamma\tau}+\frac{3}{4}\big(\pi_{\alpha\gamma}\pi_{\sigma\tau}+\pi_{\alpha\tau}\pi_{\sigma\gamma}\big)\Big)+\ \frac{3}{4}\,g_{\alpha\bar\beta}\,g_{\gamma\bar\delta}\\
				&=-\pi^\sigma_{\,.\,\bar\beta}\pi^\tau_{\,.\,\bar\delta}\,\hat S_{\alpha\sigma\gamma\tau}-\frac{3}{4}\,\pi_{\alpha\gamma}\,\pi_{\bar\beta\bar\delta}+\frac{3}{4}\,g_{\alpha\bar\delta}\,g_{\gamma\bar\beta}+\ \frac{3}{4}\,g_{\alpha\bar\beta}\,g_{\gamma\bar\delta}\\
				&=-\pi^\sigma_{\,.\,\bar\beta}\pi^\tau_{\,.\,\bar\delta}\,\hat S_{\alpha\sigma\gamma\tau}\ +\ \frac{3}{2}g\Big(R_0(e_\alpha,e_{\bar\beta})e_\gamma, e_{\bar\delta}\Big),
			\end{aligned}
		\end{equation*}
		where $\hat S$ is given by \eqref{new-dis}. If  $\mathcal K$ is as in Lemma~\ref{1-1 HKC}, we obtain that
		\begin{equation*}
			R'=-\mathcal K (\hat S)\qquad\text{and}\qquad R=-\mathcal K (\hat S)+\frac{3}{2}\,R_0.
		\end{equation*}
	\end{proof}

	We can deal with the split real form $g_{2(2)}$ of $g_2\otimes \mathbb C$ in a similar fashion. Of course, here we need to consider a different complex conjugation map; for example, a conjugation map that preserves  
	$\psi_1,\psi_2,\psi_3$, $\varphi_1, \varphi_2, \varphi_3$ and maps $e_\alpha$ to $-e_{\bar\alpha}$. If we replace the vectors $e_\alpha$ and $e_{\bar\alpha}$ in the basis \eqref{g_2_basis} of $g_2\otimes \mathbb C$ with $f_\alpha=i\,e_\alpha$ and $f_{\bar\alpha}=i\,e_{\bar\alpha}$, then the structure equations \eqref{g2_dpsi}, \eqref{g2_dvarphi} and \eqref{dealpha} take the form
	
	\begin{equation}\label{g2_dpsi_sec}
		\begin{aligned}
			d\psi^1&\ =\ -\psi^2\wedge\psi^3\ -\ \pi^\sigma_{\,.\,\bar\beta}\,(\Upsilon_1)_{\alpha\sigma}\,f^{\alpha}\wedge f^{\bar\beta}\\
			d\psi^2&\ =\ -\psi^3\wedge\psi^1\ -\ \pi^\sigma_{\,.\,\bar\beta}\,(\Upsilon_2)_{\alpha\sigma}\,f^{\alpha}\wedge f^{\bar\beta}\\
			d\psi^3&\ =\ -\psi^1\wedge\psi^2\ -\ \pi^\sigma_{\,.\,\bar\beta}\,(\Upsilon_3)_{\alpha\sigma}\,f^{\alpha}\wedge f^{\bar\beta}\\
		\end{aligned}
	\end{equation}

	\begin{equation}\label{g2_dvarphi_sec}
		\begin{aligned}
			d\varphi^1&\ =\ -\varphi^2\wedge\varphi^3\ +\ \frac{3i}{2}\Big(f^1\wedge f^{\bar 1}+f^2\wedge f^{\bar 2}+f^3\wedge f^{\bar 3}+f^4\wedge f^{\bar 4}\Big),\\
			d\varphi^2&\ =\ -\varphi^3\wedge\varphi^1\ +\ \frac{3}{2}\Big(f^1\wedge f^{3}+f^{2}\wedge f^{4}+f^{\bar 1}\wedge f^{\bar 3}+f^{\bar 2}\wedge f^{\bar 4}\Big),\\
			d\varphi^3&\ =\ -\varphi^1\wedge\varphi^2\ -\ \frac{3i}{2}\Big(f^1\wedge f^{3}+f^{2}\wedge f^{4}-f^{\bar 1}\wedge f^{\bar 3}-f^{\bar 2}\wedge f^{\bar 4}\Big),
		\end{aligned}
	\end{equation}

	\begin{equation}\label{dealpha_sec}
		\begin{aligned}
			df^\alpha&\ =\ -\pi^{\alpha\sigma}\,(\Upsilon_s)_{\sigma\beta}\,\psi^s\wedge f^\beta\ -\ \frac{i}{2}\,\varphi^1\wedge f^\alpha\ +\ \frac{1}{2}\,\pi^\alpha_{.\, \bar\beta}\,\big(\varphi^2+i\varphi^3\big)\wedge f^{\bar\beta}, \\
			df^{\bar\alpha}&\ =\ -\pi^{\bar\alpha\bar\sigma}\,(\Upsilon_s)_{\bar\sigma\bar\beta}\,\psi^s\wedge f^{\bar\beta}\ +\ \frac{i}{2}\,\varphi^1\wedge f^{\bar\alpha}\ +\ \frac{1}{2}\,\pi^{\bar\alpha}_{.\, \beta}\,\big(\varphi^2-i\varphi^3\big)\wedge f^{\beta},
		\end{aligned}
	\end{equation}

	Similarly, we may regard $G_{2(2)}$ as a principal $SO(4)$-bundle over  the symmetric space $G_{2(2)}\big/ SO(4)$. The tautological 1-form is given by $(f^\alpha,f^{\bar\alpha})$. The 1-forms $(\psi_s,\phi_s)$ represent the Levi-Civita connection. Using the same argument as in the proof of Proposition~\ref{wolf_curvature}, we can characterize the curvature tensor and the respective cubic discriminant on $G_{2(2)}\big/ SO(4)$ as follows:
	
	\begin{prop}\label{wolf_curvature} If $R=R'+\frac{Scal}{64}R_0$ is the splitting, in accordance with \eqref{gen_cur_split},  of the Riemannian curvature tensor of the symmetric space $G_{2(2)}/SO(4)$, then $R'$ is a parallel cubic discriminant.
	\end{prop}

	\section{The canonical connection} \label{sec_can_con}
	Consider an 8-dimensional quaternion-K\"ahler manifold, $(M,g,\mathcal Q)$, equipped with a cubic discriminant, $K$. Let $\nabla$ be the Levi-Civita connection of the metric. We can define a new tensor field $t$ on $M$ by the formula
	\begin{equation}\label{def_t_field}
		\begin{aligned}
			t(x,y,z)=\frac{1}{120}\sum_{a,b,c=1}^8\Big((\nabla_x K)(y,h_a,&h_b,h_c)\,K(z,h_a,h_b,h_c)\\
			& -\ (\nabla_x K)(z,h_a,h_b,h_c)\,K(y,h_a,h_b,h_c)\Big),
		\end{aligned}
	\end{equation}
	where $x$, $y$ and $z$ are vector fields and $\{h_1,\dots,h_{8}\}$ is any local $g$-orthonormal framing on $M$. 
	We shall call the new connection $\overset{c}{\nabla}$, defined by
	\begin{equation}\label{form_can_con}
		\overset{c}{\nabla}_x\,y\ =\ \nabla_x\,y\ +\ \sum_{a=1}^8 t(x,y,h_a)h_a,
	\end{equation}
	the canonical connection of $(M,g,\mathcal Q,K)$. 
	
	\begin{prop}\label{prop_Hol_can}
		The holonomy group of the canonical connection $\overset{c}{\nabla}$  is contained in $SO(4)_{ir}\subset Sp(2)Sp(1)$, i.e., we have  
		\begin{equation*}
			\overset{c}{\nabla}_x g=0,\qquad \overset{c}{\nabla}_x\mathcal Q\subset \mathcal Q, \quad\text{and}\quad\overset{c}{\nabla}_x K=0\qquad\forall x\in TM.
		\end{equation*}
		The torsion tensor,  $\overset{c}{T}$, associated with the connection $\overset{c}{\nabla}$  is defined as:
		\begin{equation*}
			\overset{c}{T}(x,y)=\overset{c}{\nabla}_xy-\overset{c}{\nabla}_yx-[x,y],\qquad x,y\in TM.
		\end{equation*}
		When we consider the action of the group $SO(4)_{ir}$ on the torsion tensor,  it splits into four components.  These components correspond to the real parts of four distinct, irreducible complex representations of $SO(4)$ (cf. Section~\ref{irr_rep_so4}) :
		\begin{equation*}
			S^9E\otimes H,\quad S^7E\otimes H,\quad S^5E\otimes H,\quad S^3E\otimes H.
		\end{equation*}
	\end{prop}
	
	Before we prove the proposition, we need to establish the following:
	\begin{lemma} 
		The tensor field $t$, defined by \eqref{def_t_field}, satisfies several key identities. These identities hold for all vector fields $x,y,z,u,v\in TM$ and $s=1,2,3$. The identities are defined with respect to any local generators $\mathcal E_s$ for the bundle $Bsp(1)_{ir}$, as in \eqref{gen_E}, and $I_s$ for the bundle $\mathcal Q$, as in \eqref{quat_id}. These identities are: 
		\begin{equation}\label{min_con_2} 
			t(x,y,z)\ =\ t(x,I_sy,I_sz);
		\end{equation}
		\begin{equation}\label{min_con_3}
			\sum_{a=1}^{8}t\big(x, h_a,\mathcal E_s\,h_a\big) \ =\ 0;
		\end{equation}
		\begin{equation}\label{min_con_4}
			\sum_{a,b=1}^{8}K(y,z,h_a,h_b)\,t(x,h_a,h_b) \ =\ -\ \frac{6}{2}\,t(x,y,z);
		\end{equation}
		\begin{equation}\label{min_con_5}
			\begin{aligned}
				(\nabla_xK)(y,z,u,v)\ =\ \sum_{a=1}^8 \Big(&t(x,y,h_a)K(h_a,z,u,v)\ +\ t(x,z,h_a)K(y,h_a,u,v)\\
				&+\ t(x,u,h_a)K(y,z,h_a,v)\ +\ t(x,v,h_a)K(y,z,u,h_a)\Big).
			\end{aligned}
		\end{equation}
	\end{lemma} 
	\begin{proof}
		For a point $p\in M$ and a given $x\in T_pM$, consider the endomorphism $H\in End(T_pM)$, 
		$$Hy=\sum_{a=1}^8t(x,y,h_a)h_a\qquad y\in T_pM.$$ 
		
		If $L$ is the covariant derivative of $K$ with respect to $x$, $L=\nabla_xK$, then, since the holonomy group of $\nabla$ preserves $\mathcal C_p\subset \mathfrak R^{Sp(n)}$, $L$ belongs to the tangent space of $\mathcal C_p$ at $K_p$. The rest follows from Lemma~\ref{tangent_space_C}.

	\end{proof}

	\begin{proof}[Proof of Proposition~\ref{prop_Hol_can}]
		The identity $t(x,y,z)=-t(x,z,y)$ implies that $\overset{c}{\nabla}_xg=0$. By \eqref{min_con_2}, $\overset{c}{\nabla}_x\mathcal Q\subset\mathcal Q$. Using \eqref{min_con_5}, we calculate 
		\begin{equation*}
			\begin{aligned}
				(\overset{c}{\nabla}_xK)(y,z,u,v)&=({\nabla}_xK)(y,z,u,v)-\sum_{a=1}^8 \Big(t(x,y,h_a)K(h_a,z,u,v)\\
				& +\ t(x,z,h_a)K(y,h_a,u,v)\ 
				+\ t(x,u,h_a)K(y,z,h_a,v)\\
				&\hspace{4cm} +\ t(x,v,h_a)K(y,z,u,h_a)\Big)\ =\ 0.
			\end{aligned}
		\end{equation*}
		Hence, the holomomy  of $\overset{c}{\nabla}$ is contained in $SO(4)_{ir}\subset Sp(2)Sp(1)$.
		
		The Lie algebra $sp(2)$ can be viewed as a representation of  $SO(4)_{ir}.$ In this context, it is equivalent to the real part of the direct sum $S^2E\oplus S^6E.$ Specifically, the  $S^2 E$ component corresponds to the subspace $sp(1)_{ir},$ while  $S^6E$ is its orthogonal complement. Therefore, by \eqref{min_con_3}, 
		the tensor $t$ takes values in a subbundle of $TM^\ast\otimes \Lambda^2(TM^\ast)$, generated by the $SO(4)$-representation 
		\begin{equation}\label{split_t}
			S^3E\otimes H\otimes S^6 E\ =\ \Big(S^9E+S^7E+S^5E+S^3E\Big)\otimes H.
		\end{equation}
		Since 
		\begin{equation*}
			g(\overset{c}{T}(x,y),w)=t(x,y,z)-t(y,x,z),\qquad \forall x,y,z\in TM,
		\end{equation*}
		the decomposition \eqref{split_t} applies also for the torsion.
	\end{proof}
	
	The next proposition establishes that the canonical connection is uniquely determined by the properties of its torsion tensor. We also prove that the underlying structure must be quaternion-K\"ahler for these specific torsion properties to hold true.
	\begin{prop}
		Let us assume we have an arbitrary 8-dimensional quaternion-Hermitian manifold  $(M,g,\mathcal Q)$ equipped with a cubic discriminant, $K$.  
		
		Now, suppose there is a linear connection $D$ on $M$ that satisfies two conditions:
		\begin{itemize}
			\item[1.] Its holonomy group is a subgroup of $SO(4)_{ir}\subset Sp(2)Sp(1).$
			\item[2.] Its torsion tensor, $T$, is related to some tensor field $d$ by the equation: 
			\begin{equation}\label{conD_con_1}
				g\big(T(x,y),z\big)=d(x,y,z)-d(y,x,z)
			\end{equation}
			for all vector fields $x$, $y$ and $z$ on $M$. The tensor field $d$ must also satisfy the following two properties for $s=1,2,3$:
			\begin{itemize}
				\item[2.1.]
				$d(x,y,z)\ =\ d(x,I_sy,I_sz).$
				\item[2.2.] $\sum_{a=1}^{8}d\big(x, h_a,\mathcal E_sh_a\big) \ =\ 0.$
			\end{itemize}
		\end{itemize}
		If all these conditions are met, then the manifold $(M,g,\mathcal Q)$ must be quaternion-K\"ahler.  Furthermore, the connection $D$ is equal to the canonical connection, $\overset{c}{\nabla}$, and the tensor field $d$ is identical to the tensor field $t$ defined in formula \eqref{def_t_field}.
	\end{prop}

	\begin{proof}
		Suppose that $D$ is a connection on $M$ with $Hol(D)\subset SO(4)_{ir}$, whose torsion $T$ satisfies \eqref{conD_con_1}.
		Consider a connection $\hat D$,
		\begin{equation*}
			\hat D_xy\overset{def}{=}D_xy-d(x,y,h_a)h_a,\qquad \forall x,y\in TM.
		\end{equation*}
		Then, the properties of $d$ imply that $Hol(\hat D)\subset Sp(2)Sp(1)$. For the torsion tensor $\hat T$ of $\hat D$, we have  
		\begin{equation*}
			\begin{aligned}
				g\big(\hat T(x,y),z\big)\ &=\ g\Big( D_xy-d(x,y,h_a)h_a- D_yx+d(y,x,h_a)h_a,z-[x,y],z\Big)\\
				&=\ g\big(T(x,y),z\big) -d(x,y,z)+d(y,x,z)\ =\ 0.
			\end{aligned}
		\end{equation*}
		Hence, $(M,g,\mathcal Q)$ is quaternion-K\"ahler and $\hat D=\nabla$ (the Levi-Civita connection). By assumption $DK=0$ and therefore, 
		\begin{equation*}
			\begin{aligned}
				(\nabla_xK)(y,z,u,v)\ =\ \sum_{a=1}^8 \Big(&d(x,y,h_a)K(h_a,z,u,v)\ +\ d(x,z,h_a)K(y,h_a,u,v)\\
				&+\ d(x,u,h_a)K(y,z,h_a,v)\ +\ d(x,v,h_a)K(y,z,u,h_a)\Big).
			\end{aligned}
		\end{equation*}
		
		For a point $p\in M$ and a given $x\in T_pM$, consider the two endomorphisms $H$ and $U$ from $End(T_pM)$, 
		$$Hy=\sum_{a=1}^8t(x,y,h_a)h_a,\qquad Uy=\sum_{a=1}^8d(x,y,h_a)h_a.$$
		By Lemma~\ref{tangent_space_C}, $H=\frac{1}{5}\Big(\frac{7}{2}U-\mathcal T_K(U)\Big)$. By assumption 2.2. of the proposition, $U\in \{\mathcal E_1,\mathcal E_2,\mathcal E_3\}^\perp$ and therefore, $\mathcal T_K(U)=-\frac{3}{2}U$ (cf. Remark~\ref{sp2_T_K_eigenspaces}). Hence, $H=U$, $d=t$ and $D=\overset{c}{\nabla}$.
	\end{proof}

	\section{Integrable cubic discriminants} Consider an 8-dimensional almost quaternion-Hermitian manifold $(M,g,\mathcal Q).$ A cubic discriminant $K$ on $M$ is called integrable, if there exists a torsion free connection $\nabla$  that preserves the structure $(g,\mathcal Q,K)$, i.e., such that $\nabla_x g=0$, $\nabla_x Q\subset Q$ and $\nabla_x K=0$ for all $x\in TM$. Clearly, if $K$ is integrable, then $(M,g,\mathcal Q)$ is quaternion-K\"ahler and the canonical connection from Section~\ref{sec_can_con} coincides with the Levi-Civita connection. The purpose of this section is to determine locally (up to a diffeomorphism) all possible cubic discriminants that are integrable. 
	
	\begin{thrm}\label{integrable_cubics} An almost quaternion-Hermitian 8-manifold $(M,g,\mathcal Q)$ admits an integrable cubic discriminant iff up to a normalizing constant it is locally isometric to one of the following three spaces:
		
		\vspace{0.2cm}
		\par (i) The flat space $\mathbb R^8$;
		\par (ii) The symmetric space $G_2\big/SO(4)$;
		\par (iii) The symmetric space $G_{2(2)}\big/SO(4)$.
	\end{thrm}
	
	\begin{proof} Suppose that $(M,g,\mathcal Q)$ is an almost quaternion-Hermitian 8-manifold endowed with an integrable cubic discriminant $K$. We shall consider the respective torsion free connection $\nabla$ as a  $so(4)_{ir}$-valued 1-form on the corresponding principal $SO(4)_{ir}$-bundle $\mathcal D(M)\rightarrow M$ (cf. \eqref{DM}).  Clearly, $D(M)$ is a subbundle of the principal $Sp(2)Sp(1)$-bundle $P(M)\rightarrow M$ determined by the underlying almost quaternion-Hermitian structure. According to the conventions introduced in Section~\ref{almostQH}, the complexification of the 8 dimensional vector space $V$ splits as $V^{\mathbb C}=W\oplus\overline{W}$. We shall work with a fixed element $\hat S\in\mathfrak S^{Sp(n)}$ that is given by \eqref{new-dis} w.r.t. a fixed $Sp(2)$-adapted basis $\{e_\alpha\}$ of $W$. The structure group $SO(4)_{ir}\subset Sp(2)Sp(1)$ of $D(M)$ may be identified with the stabilizer of $\mathcal K(\hat S)$  within $Sp(2)Sp(1)$, cf. \eqref{sp2sp1-so4}. As before, the tautological 1-form $\theta$ on $\mathcal D(M)$ takes values in $V$. The composition of $\theta$ with $e^{\alpha}\in W^\ast$ and $e^{\bar\alpha}\in\overline{W}^\ast$ produces a set of globally defined 1-forms $\theta^{\alpha},\theta^{\bar\alpha}$ on $\mathcal D(M)$.
		
		If we think of the Lie algebra of $so(4)_{ir}$ as a sum $sp(1)_{ir}\oplus sp(1)$, then the torsion free connection $\nabla$ on $\mathcal D(M)$ is represented by a pair of $End(V)$-valued 1-forms:
		\begin{equation}\label{torfree_connection}
			\Big(\psi^s\, E_s,\ \frac{1}{2}\,\varphi^s J_s \Big),
		\end{equation}
		where $\psi^1$, $\psi^2$, $\psi^3$ and  $\varphi^1$, $\varphi^2$, $\varphi^3$ are globally defined real-valued 1-forms on $\mathcal D(M)$;  $E_s$ are the matrices given by \eqref{rep_W} and $J_s$ are as in \eqref{Levi-Civita}. 
		Observe, that the action of the Lie algebra $sp(1)_{ir}\oplus sp(1)$ on $V$ is such that the pair \eqref{torfree_connection} maps any vector $\theta^\alpha\,e_\alpha + \theta^{\bar\alpha}\,e_{\bar\alpha}\in V$ to a vector $\eta^\alpha\,e_\alpha + \eta^{\bar\alpha}\,e_{\bar\alpha}\in V$, where
		\begin{equation*}
			\eta^\alpha\ =\ \pi^{\alpha\sigma}\,(\Upsilon_s)_{\sigma\beta}\,{\psi}^s\,\theta^{\beta}\ +\ \frac{i}{2}\,{\varphi}^1\,\theta^\alpha\ -\ \frac{1}{2}\,\pi^\alpha_{\,.\, \bar\beta}\,\big({\varphi}^2+i{\varphi}^3\big)\,\theta^{\bar\beta}
		\end{equation*}
		and $\eta^{\bar\alpha}=\overline{\eta^\alpha}$, where $\Upsilon_s$ are the matrices \eqref{upsilon}.  Therefore, we have the structure equations
		\begin{equation}\label{torfree_str_eq}
			\begin{aligned}
				d\theta^{\alpha}\  =\ &-\pi^{\alpha\sigma}\,(\Upsilon_s)_{\sigma\beta}\,\psi^s\wedge\theta^{\beta}\ -\ \frac{i}{2}\,\varphi^1\wedge\theta^\alpha\ +\ \frac{1}{2}\,\pi^\alpha_{\,.\, \bar\beta}\,\big(\varphi^2+i\varphi^3\big)\wedge\theta^{\bar\beta}.
			\end{aligned}
		\end{equation} 
		
		The idea of the proof is to show that, up to a free choice of a constant, we can determine explicitly the rest of the structure equations concerning the six exterior derivatives $d\psi_s$ and $d\varphi_s$ on $\mathcal D(M)$.

		Applying the exterior derivative to both sides of \eqref{torfree_str_eq}, we get
		\begin{equation}\label{torfree_str_bi}
			\begin{aligned}
				0\  =\ &-\pi^{\alpha\sigma}\,\Big(\sum_{(ijk)\footnotemark[1]}(\Upsilon_i)_{\sigma\beta}\,\big(d\,\psi^i+\psi^j\wedge\psi^k\big)\ -\ \frac{i}{2}\,\pi_{\sigma\beta}\,\big(d\,\varphi^1+\varphi^2\wedge\varphi^3\big)\Big)\wedge\theta^{\beta}\\\
				&\hspace{4.5cm} +\ \frac{1}{2}\,\pi^\alpha_{\,.\, \bar\beta}\,\Big(d\,\varphi^2+\varphi^3\wedge\varphi^1+i\,d\,\varphi^3+i\,\varphi^1\wedge\varphi^2\Big)\wedge\theta^{\bar\beta},
			\end{aligned}
		\end{equation}
		\footnotetext[1]{This indicates a summation over all positive permutations $(ijk)$ of $(123)$.}
		\noindent This equation is in fact the first Bianchi identity for the connection \eqref{torfree_connection} and we shall use it to determine the exterior derivatives $d\psi_s$ and $d\varphi_s$ (i.e., the curvature of the connection). 
		
		To begin with, observe that \eqref{torfree_str_bi} implies that there exist indexed families $ (C^s)_{\alpha\beta}=- (C^s)_{\beta\alpha}$, $ (D^s)_{\alpha\bar\beta}=-(D^s)_{\bar\beta\alpha}$, $ (F^s)_{\alpha\beta}=- (F^s)_{\beta\alpha}$ and $ (G^s)_{\alpha\bar\beta}=-(G^s)_{\bar\beta\alpha}$ of functions on  $\mathcal D(M)$, so that
		\begin{equation}\label{torfree_str_bi_sym}
			\begin{aligned}
				d\,\psi^i+\psi^j\wedge\psi^k\ & =\ \frac{1}{2}\,(C^i)_{\alpha\beta}\,\theta^\alpha\wedge\theta^\beta\ + \ (D^i)_{\alpha\bar\beta}\,\theta^\alpha\wedge\theta^{\bar\beta}\ +\ \frac{1}{2}\,(C^i)_{\bar\alpha\bar\beta}\,\theta^{\bar\alpha}\wedge\theta^{\bar\beta}\\
				d\,\varphi^i+\varphi^j\wedge\varphi^k\ &=\ \frac{1}{2}\,(F^i)_{\alpha\beta}\,\theta^\alpha\wedge\theta^\beta\ + \ (G^i)_{\alpha\bar\beta}\,\theta^\alpha\wedge\theta^{\bar\beta}\ +\ \frac{1}{2}\,(F^i)_{\bar\alpha\bar\beta}\,\theta^{\bar\alpha}\wedge\theta^{\bar\beta}.
			\end{aligned}
		\end{equation}
		Substituting \eqref{torfree_str_bi_sym} into \eqref{torfree_str_bi}, we consider the total sum of coefficients for each of the term types $\theta^{\beta}\wedge\theta^{\gamma}\wedge\theta^{\delta}$; $\theta^{\beta}\wedge\theta^{\gamma}\wedge\theta^{\bar\delta}$, $\theta^{\beta}\wedge\theta^{\bar\gamma}\wedge\theta^{\bar\delta}$ and $\theta^{\bar\beta}\wedge\theta^{\bar\gamma}\wedge\theta^{\bar\delta}$, separately, as an equation. We obtain, respectively, the following identities:
		\begin{equation}\label{cf_type_1}
			\begin{aligned}
				(\Upsilon_s)_{\alpha\beta}\,(C^s)_{\gamma\delta}+(\Upsilon_s)_{\alpha\gamma}\,&(C^s)_{\delta\beta}+(\Upsilon_s)_{\alpha\delta}\,(C^s)_{\beta\gamma} \\
				&\ =\ \frac{i}{2}\Big(\pi_{\alpha\beta}\,(F^1)_{\gamma\delta}+\pi_{\alpha\gamma}\,(F^1)_{\delta\beta}+\pi_{\alpha\delta}\,(F^1)_{\beta\gamma}\Big);
			\end{aligned}
		\end{equation}

		\begin{equation}\label{cf_type_2}
			\begin{aligned}
				(\Upsilon_s)_{\alpha\beta}\,(D^s)_{\gamma\bar\delta} - (\Upsilon_s)_{\alpha\gamma}\,&(D^s)_{\beta\bar\delta}\\
				&\ =\ \frac{1}{2}\ g_{\alpha\bar\delta}\,(F^2+iF^3)_{\beta\gamma} \ +\ \frac{i}{2}\Big(\pi_{\alpha\beta}\,(G^1)_{\gamma\bar\delta}-\pi_{\alpha\gamma}\,(G^1)_{\beta\bar\delta}\Big);
			\end{aligned}
		\end{equation}
		
		\begin{equation}\label{cf_type_3}
			\begin{aligned}
				(\Upsilon_s)_{\alpha\beta}\,(C^s)_{\bar\gamma\bar\delta}\ =\   \frac{i}{2}\,\pi_{\alpha\beta}\,(F^1)_{\bar\gamma\bar\delta}\ +\  \frac{1}{2}\,g_{\alpha\bar\delta}\,( G^2+i\,G^3)_{\beta\bar\gamma}\  -\  \frac{1}{2}\,g_{\alpha\bar\gamma}\,(G^2+i\,G^3)_{\beta\bar\delta};
			\end{aligned}
		\end{equation}

		\begin{equation}\label{cf_type_4}
			\begin{aligned}
				g_{\alpha\bar\beta}\, (F^2+i\,F^3)_{\bar\gamma\bar\delta}+g_{\alpha\bar\gamma}\, (F^2+i\,F^3)_{\bar\delta\bar\beta}+g_{\alpha\bar\delta}\, (F^2+i\,F^3)_{\bar\beta\bar\gamma}\ =\ 0.
			\end{aligned}
		\end{equation}
		
		Contracting with $\pi^{\alpha\beta}$ in \eqref{cf_type_3}, we get
		\begin{equation}\label{contr_cf_type_3}
			(F^1)_{\alpha\beta}\ =\ -\frac{1}{4}\pi^{\bar\sigma}_{\ \alpha}\,(G^3+i\,G^2)_{\bar\sigma\beta}\ +\ \frac{1}{4}\pi^{\bar\sigma}_{\ \beta}\,(G^3+i\,G^2)_{\bar\sigma\alpha}.
		\end{equation}

		Let us multiply the complex conjugate of \eqref{cf_type_3} by $\pi^{\bar\alpha}_{\ \sigma}\pi^{\bar\beta}_{\ \tau}$ and sum over the indices $\bar\alpha$, $\bar\beta$. Then, since $(\mathfrak j \Upsilon_s)_{\alpha\beta}=(\Upsilon_s)_{\alpha\beta}$, we obtain the equation
		\begin{equation}\label{cf_type_3_sec}
			\begin{aligned}
				(\Upsilon_s)_{\alpha\beta}\,(C^s)_{\gamma\delta}\ =\   -\frac{i}{2}\,\pi_{\alpha\beta}\,(F^1)_{\gamma\delta}\ -\  \frac{1}{2}\,&\pi_{\alpha\delta}\,\pi^{\bar\sigma}_{\ \beta}\,( G^2-i\,G^3)_{\bar\sigma\gamma}\\
				&\qquad  +\   \frac{1}{2}\,\pi_{\alpha\gamma}\,\pi^{\bar\sigma}_{\ \beta}\,( G^2-i\,G^3)_{\bar\sigma\delta}.
			\end{aligned}
		\end{equation}
		
		Comparing the cyclic sum over $\alpha$, $\beta$, $\gamma$ of \eqref{cf_type_3_sec} with \eqref{cf_type_1}, we get
		\begin{equation*}
			(F^1)_{\alpha\beta}\ =\ \frac{1}{2}\pi^{\bar\sigma}_{\ \alpha}\,(G^3+i\,G^2)_{\bar\sigma\beta}\ -\ \frac{1}{2}\pi^{\bar\sigma}_{\ \beta}\,(G^3+i\,G^2)_{\bar\sigma\alpha},
		\end{equation*}
		and therefore, in view of \eqref{contr_cf_type_3}, we obtain
		\begin{equation*}\label{vanFG}
			(F^1)_{\alpha\beta}=(C^s)_{\alpha\beta}=0,\qquad (G^2)_{\alpha\bar\beta}=(G^3)_{\alpha\bar\beta}=0.
		\end{equation*}
		By \eqref{cf_type_4}, we have also
		\begin{equation}\label{eqF23}
			(F^2)_{\alpha\beta}=i\,(F^3)_{\alpha\beta}.
		\end{equation}
		
		A cyclic summation over $\alpha$, $\beta$, $\gamma$ of both sides in \eqref{cf_type_2} gives
		\begin{equation*}
			\sum_{(\alpha\beta\gamma)}\Big(g_{\alpha\bar\delta}\,(F^2)_{\beta\gamma}\ +\ i\,\pi_{\alpha\beta}(G^1)_{\gamma\bar\delta} \Big)\ =\ 0.
		\end{equation*}
		Contracting the latter with $\pi^{\alpha\beta}$, we get
		\begin{equation}\label{G1_eq}
			(G^1)_{\alpha\bar\beta}\ =\ -i\,\pi^{\sigma}_{\ \bar\beta} (F^2)_{\sigma\alpha}\ +\ \frac{i}{2}\,\pi^{\sigma\tau}(F^2)_{\sigma\tau}\,g_{\alpha\bar\beta}.
		\end{equation}
		
		By substituting \eqref{G1_eq} and \eqref{eqF23} into \eqref{cf_type_2}, we end up with a linear equation  (w.r.t. the unknown quantities $(D^s)_{\alpha\bar\beta}$ and $(F^2)_{\alpha\beta}$ ):
		\begin{equation}\label{cf_type_2_new}
			\begin{aligned}
				(\Upsilon_s)_{\alpha\beta}\,(D^s)_{\gamma\bar\delta} - (\Upsilon_s)_{\alpha\gamma}\,(D^s)_{\beta\bar\delta}\ &=\ g_{\alpha\bar\delta}\,(F^2)_{\beta\gamma}\ +\ \frac{1}{2}\pi^{\sigma}_{\ \bar\delta}\,\Big(\pi_{\alpha\beta}(F^2)_{\sigma\gamma}-\pi_{\alpha\gamma}(F^2)_{\sigma\beta}\Big)\\ 
				&\quad -\ \frac{1}{4}\pi^{\sigma\tau}\,(F^2)_{\sigma\tau}\,\Big(\pi_{\alpha\beta}\,g_{\gamma\bar\delta}-\pi_{\alpha\gamma}\,g_{\beta\bar\delta}\Big).
			\end{aligned}
		\end{equation}
		The easiest way to deal with this is, probably, to solve it  using one's favorite computer algebra systems.  The result is that the solutions are
		\begin{equation*}
			(F^2)_{\alpha\beta}\ =\ h\,\pi_{\alpha\beta},\qquad (D^s)_{\alpha\bar\beta}\ =\ -\frac{2h}{3}\,(\Upsilon_s)_{\alpha\sigma}\,\pi^\sigma_{\ \bar\beta},
		\end{equation*}
		where $h$ is an arbitrary (unknown) real-valued function on $\mathcal D(M)$. 
		
		Summarizing, we have established, as a consequence of the first Bianchi identity \eqref{torfree_str_bi}, the following relations 
		
		\begin{equation}\label{str_psi_h}
			\begin{aligned}
				d\,\psi^1+\psi^2\wedge\psi^3\ & =\ -\frac{2h}{3}\,\pi^\sigma_{\,.\,\bar\beta}\,(\Upsilon_1)_{\alpha\sigma}\,\theta^{\alpha}\wedge \theta^{\bar\beta}\\
				d\,\psi^2+\psi^3\wedge\psi^1\ & =\ -\frac{2h}{3}\,\pi^\sigma_{\,.\,\bar\beta}\,(\Upsilon_2)_{\alpha\sigma}\,\theta^{\alpha}\wedge \theta^{\bar\beta}\\
				d\,\psi^3+\psi^1\wedge\psi^2\ & =\ -\frac{2h}{3}\,\pi^\sigma_{\,.\,\bar\beta}\,(\Upsilon_3)_{\alpha\sigma}\,\theta^{\alpha}\wedge \theta^{\bar\beta}
			\end{aligned}
		\end{equation}
		\begin{equation}\label{str_phi_h}
			\begin{aligned}
				d\,\varphi^1+\varphi^2\wedge\varphi^3\ &=\  -\ i\,h\,g_{\alpha\bar\beta}\,\theta^{\alpha}\wedge \theta^{\bar\beta},\\
				d\,\varphi^2+\varphi^3\wedge\varphi^1\ &=\  \ \frac{h}{2}\,\Big(\pi_{\alpha\beta}\,\theta^{\alpha}\wedge \theta^{\beta}\ +\ \pi_{\bar\alpha\bar\beta}\,\theta^{\bar\alpha}\wedge \theta^{\bar\beta}\Big),\\
				d\,\varphi^3+\varphi^1\wedge\varphi^2\ &=\  -\ \frac{ih}{2}\,\Big(\pi_{\alpha\beta}\,\theta^{\alpha}\wedge \theta^{\beta}\ -\ \pi_{\bar\alpha\bar\beta}\,\theta^{\bar\alpha}\wedge \theta^{\bar\beta}\Big).
			\end{aligned}
		\end{equation}
		
		In order to show that $h$ is actually a constant, we differentiate the first line in \eqref{str_phi_h} and then,  using  \eqref{torfree_str_eq} and \eqref{str_phi_h}, we simplify the resulting equation to 
		\begin{equation*}
			i\,g_{\alpha\bar\beta}\,dh\wedge\theta^\alpha\wedge\theta^{\bar\beta}\ =\ 0.
		\end{equation*}
		This yields $dh=0$.
	\end{proof}
	By rescaling the metric $g$, we can normalize the constant $h$ to one of three specific values:
	\begin{equation}
		h=-\frac{3}{2},\qquad h=0\qquad \text{or}\qquad h=\frac{3}{2}.
	\end{equation}
	Depending on the value of $h$, the structure equations \eqref{torfree_str_eq}, \eqref{str_psi_h} and \eqref{str_phi_h} take, respectively, the form of those of: $G_2/SO(4)$ (cf. \eqref{g2_dpsi}, \eqref{g2_dvarphi} and \eqref{dealpha}); the flat space; or $G_{2(2)}/SO(4)$ (cf. \eqref{g2_dpsi_sec}, \eqref{g2_dvarphi_sec} and \eqref{dealpha_sec}). This completes the proof.

	\section{Curvature characterization for $G_2/SO(4)$ and $G_{2(2)}/SO(4)$}
	The curvature tensor $R$ of any 8-dimensional quaternion-K\"ahler manifold $(M,g,\mathcal Q)$ decomposes as $R=R'+\frac{Scal}{64}R_0$, cf. \eqref{gen_cur_split}. If $R'$ vanishes, then depending on whether  $Scal>0$, $Scal<0$ or $Scal=0$,  the manifold is locally isometric (up to rescaling) to either the quaternion projective space $\mathbb HP^2$, the quaternion hyperbolic space $\mathbb HH^2$ or the flat space $\mathbb H^2$.   The purpose of this section is to establish a similar result (Theorem~\ref{main_thrm_third} below) where instead of assuming that $R'$ vanishes, we assume that $R'$ is pointwise proportional to a cubic discriminant. The proof is based on Theorem~\ref{integrable_cubics} from the previous section.

	\begin{thrm}\label{main_thrm_third} Suppose the curvature tensor $R$ of a connected 8-dimensional quaternion-K\"ahler manifold  $(M,g,\mathcal Q)$ decomposes as $R=R'+\frac{Scal}{64}R_0$, as discussed in \eqref{gen_cur_split}.  If there is a function $f$ on $M$ such that  $f R'$ is a cubic discriminant, then $f$ is constant. In this case, depending on the sign of the scalar curvature ($Scal$), the manifold $M$ is locally isometric (up to rescaling) to either $G_2/SO(4)$ (if $Scal>0$) or $G_{2(2)}/SO(4)$ (if $Scal<0$).
	\end{thrm}

	Before starting the proof of the theorem, let us make a few remarks: 
	\begin{rmrk}
		Note that the theorem does not require a non-zero scalar curvature. This is important because hyper-K\"ahler manifolds, which are a special case of quaternion-K\"ahler manifolds, have zero scalar curvature. Therefore, the theorem implies that a hyper-K\"ahler manifold's curvature can never be pointwise proportional to a cubic discriminant, except when the manifold is the flat space.
	\end{rmrk}
	
	Notice that the condition for $f R'$ to be a cubic discriminant on a quaternion-K\"ahler manifold implies that $f$ is a nowhere vanishing function. Additionally, according to Theorem~\ref{sp2-characterization},   this condition holds if and only if both of the following are true:
	\vspace{0.1cm}
	
	(I) $\ $ $(2f\mathcal T_{R'}-7\,\text{Id})(2f\mathcal T_{R'}+3\,\text{Id})\ =\ 0.$
	
	\vspace{0.1cm}
	
	(II) $\ $ For any sections $A$ and $B$ of $Bsp(2)\subset End(TM)$ (cf. \eqref{Bsp2}):
	$$f\big[\mathcal T_{R'}A,\mathcal T_{R'}B\big]-f\mathcal T_{R'}\big[\mathcal T_{R'}A, B\big]\ =\ \frac{3}{2}\,\Big(\mathcal T_{R'}\big[A,B\big]-\big[A,\mathcal T_{R'}B\big]\Big).$$

	In our proof, we shall use the notation from Section~\ref{almostQH}. We shall consider the principal $Sp(2)Sp(1)$-bundle $P(M)\rightarrow M$ with its canonical 1-forms, $\theta^{\alpha}$ and $\theta^{\bar\alpha}$. The Levi-Civita connection is defined by the 1-forms $\phi^s$ and  $\Gamma_{\alpha\beta}$ on $P(M)$.  These forms satisfy the following structure equations (cf. \eqref{str_eq_LC}):
	\begin{equation}\label{first_str_QK}
		d\theta^{\alpha}\ =\ -\pi^{\alpha\sigma}\,\Gamma_{\sigma\beta}\wedge\theta^{\beta}\ -\ \frac{i}{2}\,\phi^1\wedge\theta^\alpha\ +\ \frac{1}{2}\,\pi^\alpha_{\,.\, \bar\beta}\,\big(\phi^2+i\phi^3\big)\wedge\theta^{\bar\beta}.
	\end{equation} 
	
	The integrability condition $d(d\theta^{\alpha})=0$ (which is essentially  the first Binachi identity for the Levi-Civita connection), obtained by differentiating \eqref{first_str_QK}, yields the existence of a real-valued function $C$ and a $\mathfrak S^{Sp(2)}$-valued function $\L$ (with $ \L_{\alpha\beta\gamma\delta}=\L(e_\alpha,e_\beta,e_\gamma,e_\delta)$), so that
	\begin{equation}\label{first_bi_res}
		\begin{aligned}
			d\Gamma_{\alpha\beta}\ &=\ -\pi^{\sigma\tau}\Gamma_{\alpha\sigma}\wedge\Gamma_{\tau\beta}\ +\ \Big(\pi^{\sigma}_{\ \bar\delta}\, \L_{\alpha\beta\gamma\sigma}-Cg_{\alpha\bar\delta}\,\pi_{\beta\gamma}-Cg_{\beta\bar\delta}\,\pi_{\alpha\gamma}\Big)\,\theta^{\gamma}\wedge\theta^{\bar\delta}\\
			d\phi^1\ &=\ -\phi^2\wedge\phi^3\ -\ 2iC\,g_{\alpha\bar\beta}\,\theta^{\alpha}\wedge\theta^{\bar\beta}\\
			d\phi^2\ &=\ -\phi^3\wedge\phi^1\ -\ C\,\pi_{\alpha\beta}\,\theta^{\alpha}\wedge\theta^{\beta}\ -\ C\,\pi_{\bar\alpha\bar\beta}\,\theta^{\bar\alpha}\wedge\theta^{\bar\beta}\\
			d\phi^3\ &=\ -\phi^1\wedge\phi^2\ +\ iC\,\pi_{\alpha\beta}\,\theta^{\alpha}\wedge\theta^{\beta}\ -\ iC\,\pi_{\bar\alpha\bar\beta}\,\theta^{\bar\alpha}\wedge\theta^{\bar\beta}.
		\end{aligned}
	\end{equation}
	One can easily recognize that $C$ and $\L$ represent the two components of the curvature in the decomposition $R=R'+\frac{Scal}{64}R_0$,
	\begin{equation*}
		Scal=\frac{128\,C}{3},\qquad  R'=-\mathcal K(\L),
	\end{equation*}
	where $\mathcal K$ is defined in Lemma~\ref{1-1 HKC}. We have that  $f R'$ is a cubic discriminant on $M$ iff the  $\mathfrak S^{Sp(2)}$-valued function $S$ (with $S_{\alpha\beta\gamma\delta}=S(e_\alpha,e_\beta,e_\gamma,e_\delta)$),
	\begin{equation}\label{cubic_S}
		S\ {=}\ -f\L,
	\end{equation} 
	satisfies the relations (I) and (II) of Corollary~\ref{characterization_in_coordinates}.

	\begin{proof}[Proof of the Theorem~\ref{main_thrm_third}] Let us begin by differentiating equations \eqref{first_bi_res}.  A straightforward calculation shows that the integrability conditions, $d(d\Gamma_{\alpha\beta})=0$ and $d(d\phi^s)=0$, are  equivalent to the relations $dC=0$ (meaning $C$ is constant) and the following equation:
		\begin{equation}\label{sec_B_solved}
			\begin{aligned}
				d\L_{\alpha\beta\gamma\delta}\ -\  \pi^{\sigma\tau}\Big(\Gamma_{\tau\alpha}\L_{\sigma\beta\gamma\delta}\ &+\ \Gamma_{\tau\beta}\L_{\alpha\sigma\gamma\delta}\ +\ \Gamma_{\tau\gamma}\L_{\alpha\beta\sigma\delta}\ +\ \Gamma_{\tau\delta}\L_{\alpha\beta\gamma\sigma}\Big)\\
				&\hspace{2cm}=\ X_{\alpha\beta\gamma\delta\mu}\,\theta^{\mu}\ -\ \pi^{\sigma}_{\ \bar\mu}\,(\mathfrak jX)_{\alpha\beta\gamma\delta\sigma}\,\theta^{\bar\mu}.
			\end{aligned}
		\end{equation}
		These integrability conditions essentially provide the second Bianchi identity for the Levi-Civita connection. The family of functions $X_{\alpha\beta\gamma\delta\mu}$  on $P(M)$ represents the covariant derivative of the curvature on  $M$  and depends symmetrically on all indices: $\alpha$, $\beta$, $\gamma$, $\delta$ and $\mu$.

		Let us assume now that $fR'$ is indeed a cubic discriminant on $M$ and define 
		$S$ by~\eqref{cubic_S}. At any point $u\in P(M)$, the value $S(u)$ is a tensor on $V=\mathbb R^8$ that lies within the $Sp(2)$-orbit $\mathcal K^{-1}(\mathcal C)$. Using the same reasoning as for equation \eqref{LinTK_U}, we can characterize the tangent space to this orbit at $S(u)$. The tangent space $T_{S(u)}\mathcal K^{-1}(\mathcal C)$ is given by:
		\begin{equation*}
			T_{S(u)}\mathcal K^{-1}(\mathcal C)\ =\ \Big\{\big[Z,S(u)\big]_{\alpha\beta\gamma\delta}\ :\  Z\in sp(2)\Big\}.
		\end{equation*}
		Here, we have used the Lie bracket formula:
		\begin{equation*}
			[Z,S]_{\alpha\beta\gamma\delta}\ =\ -\pi^{\sigma\tau}\Big(Z_{\tau\alpha}S_{\sigma\beta\gamma\delta}\ +\ Z_{\tau\beta}S_{\alpha\sigma\gamma\delta}\ +\ Z_{\tau\gamma}S_{\alpha\beta\sigma\delta}\ +\ Z_{\tau\delta}S_{\alpha\beta\gamma\sigma}\Big)
		\end{equation*}
		for any $Z=\{Z_{\alpha\beta}\}\in sp(2)$ with the identification \eqref{Spn_X}. Therefore, there exist 1-forms:
		\begin{equation}\label{def_Z3}
			Z_{\alpha\beta}\ =\ Z_{\alpha\beta\gamma}\,\theta^{\gamma}\ -\ \pi^{\sigma}_{\ \bar\gamma}(\mathfrak j Z)_{\alpha\beta\mu}\,\theta^{\bar\gamma},
		\end{equation}
		where the family of functions $Z_{\alpha\beta\gamma}$ on $P(M)$ is symmetric in indices  $\alpha$ and $\beta$ (but not necessarily totally symmetric in all three indices). These 1-forms are defined such that they satisfy the following relation:
		\begin{equation*}
			dS_{\alpha\beta\gamma\delta}\ -\  \pi^{\sigma\tau}\Big(\Gamma_{\tau\alpha}S_{\sigma\beta\gamma\delta}\ +\ \Gamma_{\tau\beta}S_{\alpha\sigma\gamma\delta}\ +\ \Gamma_{\tau\gamma}S_{\alpha\beta\sigma\delta}\ +\ \Gamma_{\tau\delta}S_{\alpha\beta\gamma\sigma}\Big)\ =\ [Z,S]_{\alpha\beta\gamma\delta}.
		\end{equation*}
		Then, by \eqref{cubic_S} and \eqref{sec_B_solved}, we have
		\begin{equation*}
			\frac{1}{f^2}\,S_{\alpha\beta\gamma\delta}\,df \ -\ \frac{1}{f}\,[Z,S]_{\alpha\beta\gamma\delta}
			\ =\ X_{\alpha\beta\gamma\delta\mu}\,\theta^{\mu}\ -\ \pi^{\sigma}_{\ \bar\mu}\,(\mathfrak jX)_{\alpha\beta\gamma\delta\sigma}\,\theta^{\bar\mu},
		\end{equation*}
		which in view of \eqref{def_Z3} yields the equation
		\begin{equation}\label{main_eq_SZX}
			\begin{aligned}
				\frac{1}{f}\,\pi^{\sigma\tau}\Big(Z_{\tau\alpha\mu}S_{\sigma\beta\gamma\delta}\ +\ Z_{\tau\beta\mu}S_{\alpha\sigma\gamma\delta}\ &+\ Z_{\tau\gamma\mu}S_{\alpha\beta\sigma\delta}\ +\ Z_{\tau\delta\mu}S_{\alpha\beta\gamma\sigma}\Big)\\ 
				&\qquad\qquad+\ \frac{1}{f^2}\,S_{\alpha\beta\gamma\delta}\,f_{\mu}\ =\  X_{\alpha\beta\gamma\delta\mu},
			\end{aligned}
		\end{equation}
		where $df\ =\ f_\mu\,\theta^{\mu}\ +\ f_{\bar\mu}\,\theta^{\bar\mu}.$ 
		
		Since $\mathcal K(S)$ is a cubic discriminant, for each $p\in M$, there exists $u\in P(M)$ that projects to $p$ such that $S_{\alpha\beta\gamma\delta}(u)\ =\ {\hat S}_{\alpha\beta\gamma\delta}$, where $\hat S$ is given by \eqref{new-dis}. We shall fix such a $u\in P(M)$ and consider  equation \eqref{main_eq_SZX} only at this point. Next, we multiply both sides of \eqref{main_eq_SZX} by 
		$$\pi^{\beta\kappa}(\Upsilon_s)_{\kappa\epsilon}\,(\Upsilon_s)_{\xi\zeta}\,\pi^{\xi\gamma}\,\pi^{\zeta\delta}$$
		(and sum over repeated indices), with  $\Upsilon_s$ as given in \eqref{upsilon}. Using Lemma \ref{upsilon_lemma}, we get 
		\begin{equation*} 
			\begin{aligned}
				\pi^{\beta\kappa}(\Upsilon_s)_{\kappa\epsilon}\,&(\Upsilon_s)_{\xi\zeta}\,\pi^{\xi\gamma}\,\pi^{\zeta\delta}\,X_{\alpha\beta\gamma\delta\mu}\ =\ \frac{105}{8f}\,Z_{\epsilon\alpha\mu}+ \frac{105}{8f^2}\,\pi_{\alpha\epsilon}\,f_{\mu}\\
				&+\ \frac{7}{2f}\,\pi^{\sigma\tau}\,\pi^{\beta\kappa}\,Z_{\tau\beta\mu}\Big(S_{\alpha\sigma\kappa\epsilon}+\frac{3}{4}\pi_{\alpha\kappa}\pi_{\sigma\epsilon}
				+\frac{3}{4}\pi_{\alpha\epsilon}\pi_{\sigma\kappa}\Big)\\
				&+\ \frac{1}{f}\,\pi^{\sigma\tau}\,\pi^{\beta\kappa}\,\pi^{\xi\gamma}\,\pi^{\zeta\delta}\,Z_{\tau\gamma\mu}\,S_{\alpha\beta\sigma\delta}\Big(S_{\kappa\epsilon\xi\zeta}+\frac{3}{4}\pi_{\kappa\xi}\pi_{\epsilon\zeta}
				+\frac{3}{4}\pi_{\kappa\zeta}\pi_{\epsilon\xi}\Big)\\
				&+\ \frac{1}{f}\,\pi^{\sigma\tau}\,\pi^{\beta\kappa}\,\pi^{\xi\gamma}\,\pi^{\zeta\delta}\,Z_{\tau\delta\mu}\,S_{\alpha\beta\gamma\sigma}\Big(S_{\kappa\epsilon\xi\zeta}+\frac{3}{4}\pi_{\kappa\xi}\pi_{\epsilon\zeta}
				+\frac{3}{4}\pi_{\kappa\zeta}\pi_{\epsilon\xi}\Big)
			\end{aligned}
		\end{equation*}
		Applying Corollary~\ref{characterization_in_coordinates}--(I) to the right-hand side of the equation above, a short calculation yields:
		\begin{equation} \label{eq_X5_trace}
			\begin{aligned}
				\pi^{\beta\kappa}(\Upsilon_s)_{\kappa\epsilon}\,&(\Upsilon_s)_{\xi\zeta}\,\pi^{\xi\gamma}\,\pi^{\zeta\delta}\,X_{\alpha\beta\gamma\delta\mu}\\ 
				&=\ \frac{21}{f}\,Z_{\epsilon\alpha\mu}+ \frac{105}{8f^2}\,\pi_{\alpha\epsilon}\,f_{\mu}
				\ -\ \frac{6}{f}\,\pi^{\sigma\tau}\,\pi^{\nu\kappa}\,Z_{\sigma\nu\mu}\, S_{\alpha\epsilon\sigma\kappa}
			\end{aligned}
		\end{equation}
		By antisymmetrizing equation \eqref{eq_X5_trace} with respect to the free indices  $\alpha$, $\epsilon$ and $\mu$, we find that:
		\begin{equation*}
			\pi_{\alpha\epsilon}\,f_\mu\ +\ \pi_{\epsilon\mu}\,f_\alpha\ +\ \pi_{\mu\alpha}\,f_\epsilon = 0.
		\end{equation*}
		This implies that $df=0$, so $f$ must be a constant. Consequently, antisymmetrizing \eqref{eq_X5_trace} with respect to the index pair  $\alpha$, $\mu$ gives the formula:
		\begin{equation}\label{eq_ZS_anti_ab}
			Z_{\alpha\beta\gamma}-Z_{\alpha\gamma\beta}\ =\ \frac{2}{7}\,\pi^{\sigma\tau}\pi^{\mu\nu}\,\Big(Z_{\sigma\mu\gamma}\,S_{\alpha\beta\tau\nu}-Z_{\sigma\mu\beta}\,S_{\alpha\gamma\tau\nu}\Big)
		\end{equation} 
		
		Next, consider the quantities $(v^1)_\alpha$, $(v^2)_\alpha$, $(v^3)_\alpha$ and $W_{\alpha\beta\gamma}$,  
		\begin{equation}\label{def_vs_W}
			(v^s)_\gamma\ \overset{def}{=}\ \frac{1}{5}\,(\Upsilon_s)_{\sigma\tau}\,\pi^{\sigma\alpha}\,\pi^{\tau\beta}\,Z_{\alpha\beta\gamma},\qquad  W_{\alpha\beta\gamma}\ \overset{def}{=}\ Z_{\alpha\beta\gamma}\ -\ \sum_s(\Upsilon_s)_{\alpha\beta}(v^s)_\gamma.
		\end{equation}
		Substituting these into \eqref{eq_ZS_anti_ab} and using Lemma~\ref{upsilon_lemma}--(5), we calculate that:
		\begin{equation*}
			\begin{aligned}
				W_{\alpha\beta\gamma}\ &-\ W_{\alpha\gamma\beta}\ +\ \sum_s(\Upsilon_s)_{\alpha\beta}(v^s)_\gamma\ -\ \sum_s(\Upsilon_s)_{\alpha\gamma}(v^s)_\beta\\
				&=\  \frac{2}{7}\,\pi^{\sigma\tau}\pi^{\mu\nu}\,Z_{\sigma\mu\gamma}\,\Bigg(\sum_s(\Upsilon_s)_{\alpha\beta}\,(\Upsilon_s)_{\tau\nu}-\frac{3}{4}\,\pi_{\alpha\tau}\,\pi_{\beta\nu}-\frac{3}{4}\,\pi_{\alpha\nu}\,\pi_{\beta\tau}\Bigg)\\
				&\hspace{1.5cm}-\ \frac{2}{7}\,\pi^{\sigma\tau}\pi^{\mu\nu}\,Z_{\sigma\mu\beta}\,\Bigg(\sum_s(\Upsilon_s)_{\alpha\gamma}\,(\Upsilon_s)_{\tau\nu}-\frac{3}{4}\,\pi_{\alpha\tau}\,\pi_{\gamma\nu}-\frac{3}{4}\,\pi_{\alpha\nu}\,\pi_{\gamma\tau}\Bigg)\\
				&=\ \frac{10}{7}\,\sum_s(\Upsilon_s)_{\alpha\beta}(v^s)_\gamma\ -\ \frac{10}{7}\,\sum_s(\Upsilon_s)_{\alpha\gamma}(v^s)_\beta\ +\ \frac{3}{7}\,Z_{\alpha\beta\gamma}\ -\ \frac{3}{7}\,Z_{\alpha\gamma\beta}.
			\end{aligned}
		\end{equation*}
		This simplifies to
		$
		W_{\alpha\beta\gamma}\ -\ W_{\alpha\gamma\beta}\ =\ 0,
		$
		which means that $W_{\alpha\beta\gamma}$ is totally symmetric. It also has the following property:
		\begin{equation}\label{W_notrace}
			(\Upsilon_s)_{\sigma\tau}\,\pi^{\sigma\alpha}\,\pi^{\tau\beta}\,W_{\alpha\beta\gamma}\ =\ 0,\qquad s=1,2,3.
		\end{equation}
		
		Multiplying both sides of equation \eqref{main_eq_SZX} by $\pi^{\alpha\mu}$ and then using \eqref{def_vs_W}, we get:
		\begin{equation*}
			\begin{aligned}
				0\ &=\ \pi^{\alpha\mu}\,\pi^{\sigma\tau}\Big(Z_{\tau\alpha\mu}S_{\sigma\beta\gamma\delta}\ +\ Z_{\tau\beta\mu}S_{\alpha\sigma\gamma\delta}\ +\ Z_{\tau\gamma\mu}S_{\alpha\beta\sigma\delta}\ +\ Z_{\tau\delta\mu}S_{\alpha\beta\gamma\sigma}\Big)\\
				&=\ \pi^{\alpha\mu}\,\pi^{\sigma\tau}\Big(W_{\tau\beta\mu}S_{\alpha\sigma\gamma\delta}\ +\ W_{\tau\gamma\mu}S_{\alpha\beta\sigma\delta}\ +\ W_{\tau\delta\mu}S_{\alpha\beta\gamma\sigma}\Big)\\
				&\qquad+\ \pi^{\alpha\mu}\,\pi^{\sigma\tau}\Big((\Upsilon_s)_{\tau\alpha}S_{\sigma\beta\gamma\delta}\ +\ (\Upsilon_s)_{\tau\beta}S_{\alpha\sigma\gamma\delta}\ +\ (\Upsilon_s)_{\tau\gamma}S_{\alpha\beta\sigma\delta}\\
				&\hspace{9cm} +\ (\Upsilon_s)_{\tau\delta}S_{\alpha\beta\gamma\sigma}\Big)\,(v^s)_\mu.
			\end{aligned}
		\end{equation*}
		
		By Lemma~\ref{upsilon_lemma}--(4), the last term in the previous equation is zero. We can use formula (5) from the same lemma and \eqref{W_notrace} to calculate:
		\begin{equation*}
			\begin{aligned}
				0\ &=\ \pi^{\alpha\mu}\,\pi^{\sigma\tau}\Big(W_{\tau\beta\mu}S_{\alpha\sigma\gamma\delta}\ +\ W_{\tau\gamma\mu}S_{\alpha\beta\sigma\delta}\ +\ W_{\tau\delta\mu}S_{\alpha\beta\gamma\sigma}\Big)\\
				&=\ -\frac{3}{4}\,\pi^{\alpha\mu}\,\pi^{\sigma\tau}\Big(W_{\tau\beta\mu}(\pi_{\alpha\gamma}\pi_{\sigma\delta}+\pi_{\alpha\delta}\pi_{\sigma\gamma})\ +\ W_{\tau\gamma\mu}(\pi_{\alpha\beta}\pi_{\sigma\delta}+\pi_{\alpha\delta}\pi_{\sigma\beta})\\
				&\hspace{6cm} +\ W_{\tau\delta\mu}(\pi_{\alpha\beta}\pi_{\sigma\gamma}+\pi_{\alpha\gamma}\pi_{\sigma\beta})\Big)
				\ =\ -\frac{9}{2}\,W_{\beta\gamma\delta}.
			\end{aligned}
		\end{equation*}
		This implies that $W_{\beta\gamma\delta}=0,$ which gives us $Z_{\alpha\beta\gamma}=\sum_s(\Upsilon_s)_{\alpha\beta}(v^s)_\gamma$. Substituting this back into \eqref{main_eq_SZX} and using  Lemma~\ref{upsilon_lemma} again, we find that $X_{\alpha\beta\gamma\delta\mu}=0$. This means that both the curvature $\L_{\alpha\beta\gamma\delta}$ (see \eqref{sec_B_solved}) and the cubic discriminant $S_{\alpha\beta\gamma\delta}$  are parallel with respect to the Levi-Civita connection. The final result then follows directly from  Theorem~\ref{integrable_cubics}.
	\end{proof}



\begin{thebibliography}{99}
		\bibitem{ABF} Agricola, I., Becker-Bender, J., Friedrich, T.,  {\it On the topology and the geometry of $SO(3)$-manifolds}, Ann. Global
		Anal. Geom. {\bf 40} (1), (2011), 67-84.
		
		\bibitem{Ale} Alekseevskii, D. V., {\it Classification of quaternionic spaces with
			a transitive solvable group of motions}, Math. USSR, Izv. 9
		(1975), 297--339.
		
		\bibitem{Ale2} Alekseevskii, D. V., {\it Riemannian spaces with exceptional holonomy groups}, Functional Analysis and
		Applications 2 (1968), 97--105.
		
		\bibitem{Bes} Besse, A., {\it Einstein manifolds}, Ergebnisse der Mathematik,
		vol. 3, Springer-Verlag, Berlin, 1987.
		
		
		\bibitem{BN} Bobi\'enski, M.,  Nurowski, P.,  {\it Irreducible $SO(3)$-geometries in dimension five}, J. Reine Angew. Math.
		605 (2007), 51-93.
		
		
		
		\bibitem{CF}  Chiossi, S., Fino, A., {\it Nearly integrable $SO(3)$ structures on 5-dimensional Lie groups},  J. Lie Theory
		17 (2007), 539-562
		
		
		\bibitem{CMS}  Conti, D.,  Madsen, T.B., Salamon, S., {\it Quaternionic geometry in dimension eight},  Geometry and Physics: a Festschrift
		in Honour of Nigel Hitchin, Oxford University Press, 2018
		
		
		\bibitem{Cor}  Cort\'es., V., {\it Alekseevskian spaces.}, Differ. Geom. Appl., 6(2), (1996), 129--168.
		
		\bibitem{FH} Fulton, W. and Harris, J.,  {\it Representation theory}, volume 129 of Graduate Texts in Mathematics, Springer-Verlag, New York, 1991, A first course, Reading in Mathematics.
		
		\bibitem{Fr} Friedrich, Th., {\it Dirac Operators in Riemannian Geometry}, American Mathematical Society, ISBN 978-0-8218-2055-1.
		
		\bibitem{H} Helgason, S., {\it Differential Geometry, Lie Groups, and Symmetric Spaces}, Academic Press, Inc., 1978.
		
		
		\bibitem{LeB} LeBrun, C., {\it Quaternionic-K\"ahler manifolds and conformal geometry}, Math. Ann. {\bf 284} (1989), 353--376.
		
		
		
		
		
		\bibitem{Sal} Salamon. S. M.,  {\it Quaternionic  K\"ahler manifolds}, Invent. Math. {\bf 67}, (1982), 143--171.
		
		
		\bibitem{Sal2} Salamon. S. M.,  {\it  Riemannian Geometry and Holonomy Groups}, Pitman Research Notes in Mathematics, vol. 201, Longman 1989.
		
		
		\bibitem{Sw} Swann, A.,  {\it Aspects symplectiques de la g\'eom\'etrie quaternionique},
		C. R. Acad. Sci. Paris S\'er. I Math., 308 (7), (1989), 225-228.
		
		
		\bibitem{W} Wolf, J. A.,  {\it Complex homogeneous contact manifolds and quaternionic symmetric
			spaces},
		J. Math. Mech., {\bf 14} (6), (1965), 1033--1047.
		
	\end{thebibliography}
\end{document}